\documentclass[10pt,conference]{IEEEtran}

\usepackage[none]{hyphenat}

\usepackage{cite}

\ifCLASSINFOpdf
   \usepackage[pdftex]{graphicx}
   \graphicspath{{figs/}}
   \DeclareGraphicsExtensions{.pdf,.jpeg,.png}
\else

   \usepackage[dvips]{graphicx}

   \graphicspath{{../figs/}}
 
   \DeclareGraphicsExtensions{.eps}
\fi

\usepackage[cmex10]{amsmath}

\usepackage{amsthm}

\usepackage{algorithmic}

\usepackage{array}

\ifCLASSOPTIONcompsoc
  \usepackage[caption=false,font=normalsize,labelfont=sf,textfont=sf]{subfig}
\else
  \usepackage[caption=false,font=footnotesize]{subfig}
\fi

\usepackage{url}

\usepackage{xcolor}
\usepackage{comment}
\usepackage{amssymb} 
\usepackage[ruled,vlined,english]{algorithm2e}

\usepackage{hyperref}

\newif\iffinal
\finaltrue

\usepackage[switch]{lineno}

\begin{document}
\sloppy 

\title{Two-Phase Optimization for PINN Training}

\author{López, D. M.}

\maketitle

\begin{abstract}
This work presents an algorithm for training Neural Networks where the loss function can be decomposed into two non-negative terms to be minimized. The proposed method is an adaptation of Inexact Restoration algorithms, constituting a two-phase method that imposes descent conditions. Some performance tests are carried out in PINN training.
\end{abstract}

\IEEEpeerreviewmaketitle

\section{Introduction}

Neural Networks are increasingly being used for various purposes. It is natural, therefore, that the loss functions used in training are constantly evolving. In particular, more and more problems are emerging in which the network is trained to model two or more different characteristics or objectives, each represented by distinct terms in the loss functions. A notable example of this is multitask learning \cite{Caruana1998}, widely employed in natural language processing problems \cite{collobert2008unified}, speech recognition \cite{huang2013cross}, and computer vision \cite{kokkinos2017ubernet}, such as in semantic segmentation \cite{jadon2020survey} \cite{azad2023loss}.

Another field where this happens is in Neural Networks which estimate solutions of differential problems. These are known as Physics Informed Neural Networks (PINN). Here, the solution (network output) is expected to satisfy not only the associated Partial Differential Equation (PDE) but also the initial and boundary conditions. Traditionally, a loss term is associated with each of these objectives. Thus, the network's goal of estimating solutions to differential problems is generally associated with two or three terms to be minimized.

A natural solution to these problems is considering the sum of the involved terms as the loss function. This is because, by minimizing the sum of non-negative functions, it is implicitly expected to minimize each one of them.

The most popular methods for training networks with these characteristics consider weighted sums as the loss function, using traditional optimization methods to minimize it. Most of these methods do not impose descent conditions. By not imposing descent conditions for the acceptance, or not, of a step in the optimization algorithm, we risk accepting new parameters for the network that do not improve its performance in relation to the previous ones. Thus, it is common that during the training of PINN using traditional methods, the loss functions suffer an increase in some iterations.

Our proposed method imposes descent conditions, improving the effectiveness of the training process. Additionally, it aims to simultaneously minimize two functions, enhancing the training of PINN and multiobjective networks.

\section{Related Work}
\subsection{Network Training}

In problems where the objective of the Neural Network (NN) can be represented by two or more terms that must be minimized, it is common to consider the weighted sum of these terms as the loss function. However, determining the ideal weights constitutes a new problem, often solved by training multiple models with different weight choices and selecting the best according to some criterion. However, this process is not efficient. An approximation to solve this problem is to train the model only once, considering a distribution of loss functions instead of training multiple times with a single loss function \cite{dosovitskiy2019you}.

A wide range of networks, whose objectives are described by two or more terms, are Physics Informed Neural Networks (PINN), initially proposed in \cite{raissi2019physics}. In these networks, the choice of loss function weights can be optimally made, as studied in \cite{van2022optimally}. However, this method requires a broad qualitative knowledge of the solutions, which is only sometimes possible or convenient. It is also possible to use a heuristic to make this choice \cite{van2022optimally}. For the case of PINN, more specific studies on weight choice have been conducted (see \cite{wang2022respecting}).

Once the loss function is well defined, an optimization algorithm is used for training the NN, with Adam \cite{kingma2014adam} being one of the most commonly used in PINN training.

On the other hand, the weighted sum of loss functions is not limited to PINN. In multitask learning, specific methods are also used to determine a good choice of weights \cite{kendall2018multi}.

Existing weight determination strategies are often tailored to specific networks and problems. In light of this, our contribution is a general method that tackles this limitation by minimizing two different loss functions in the same algorithm, thereby enabling the training of any network with this characteristics. 

\subsection{Inexact Restoration}

Training a Neural Network is essentially an optimization process in which we minimize the loss function and adjust the network's weights. For this purpose, any suitable optimization algorithm can be used. However, many popular methods for training Neural Networks do not impose any condition to accept or reject the step taken. This strategy of imposing descent conditions is widely used in more elaborate optimization algorithms \cite{zhang2009nonlinear, malik2020new, andrei2013another, moreno2023globalizaccao, santos2022globalizaccao, bueno2015flexible, grapiglia2017worst, martinez2000inexact, fischer2010new} and allows demonstrating convergence to local or global minima. This also intends to guarantee that the behavior of the function decreases over the optimization process. However, these methods are not the most popular in machine learning.

A broad spectrum of methods allows imposing descent conditions throughout the optimization process. However, we will focus on a set of methods that, in addition to imposing descent conditions, involve minimizing two functions in their development. We believe that by not always using the gradient information of the same function as the descent direction, these methods have less tendency to converge to local minima of the sum function. We are referring to the Inexact Restoration methods \cite{martinez2000inexact, martinez2001inexact, fischer2010new, bueno2020complexity}.

Originally, Inexact Restoration was conceived to solve constrained  optimization problems, such as \eqref{eq:optrest}:

\begin{equation} \label{eq:optrest}
	\begin{array}{cl} 
		\mbox{Minimize}  &  f(x) \\ 
		\mbox{Subject to}  &   \|h(x)\|= 0, \\
		& x \in \Omega ,
	\end{array}
\end{equation}
where $\Omega$ is a box constraint, limiting the parameter $x$ from tending to infinity. A point $x$ is said to be feasible when, besides satisfying $x \in \Omega$, it satisfies $\|h(x)\|=0$. This latter function is used to measure the feasibility of a point, so we say, for example, that $x$ is more feasible than $y$ when $\|h(x)\|\leq \|h(y)\|$. In the optimization context, the function $f$ to be minimized is called the objective function.

According to \cite{bueno2015flexible}, the main characteristics of Inexact Restoration methods can be synthesized with the following two steps:

\begin{itemize}
    \item An arbitrary method calculates a point that is sufficiently more feasible given the current iterate. This step is called the restoration phase.
    \item In the optimization phase, a function associated with the objective function is minimized, and the result is compared with the current iterate regarding feasibility and optimality. If the point is accepted, it will be set as the new iterate; otherwise, a new point is found in a region closer to the point obtained in the restoration phase.
\end{itemize}

This is a two-phase method in which the strategy aims to minimize $f(x)$, allowing points where $h(x)\neq 0$ in the optimization process, but expecting the convergence to a feasible point ($h(x^\star) = 0$). It is known that the minimum of $\|h(x)\|$ is achieved when this function is zero. If this is not possible, we will be faced with an infeasible problem, one that has no solution.

Although not designed with this view, Inexact Restoration methods essentially minimize two functions jointly. Therefore, in this work, we will present an adaptation, strongly inspired by Inexact Restoration methods, to be used in training networks with multiple (two) loss functions.

\section{Proposed Technique}

In this work proposal, we aim to minimize two functions, $L_1$ and $L_2$. Drawing an analogy with the original Inexact Restoration (IR) methods, the function $L_2$ will assume the role of $\|h(x)\|$ within the algorithm, and $L_1$ will represent $f(x)$. Thus, following the spirit of IR, the method will consist of two phases:

\begin{itemize}
    \item First phase: Given an iterate $x$, search, through some method, for a new iterate $y$ such that $L_2(y)$ is sufficiently smaller than $L_2(x)$, or equal if the minimum has been reached.
    \item Second phase: Minimize the function $L_1$, but the new step will only be accepted when criteria aimed at measuring the quality of the new iterate in terms of not only the decrease of $L_1$ but also the non-increase of $L_2$ are satisfied.
\end{itemize}

In this second phase, the descent criterion will be associated with the decrease of the penalty function:
\begin{equation}
    \Phi(x,\theta)=\theta L_1(x) + (1-\theta)L_2(x)
\end{equation}
proposed in \cite{martinez2000inexact}, where the penalty parameter $\theta$ is a real value in $(0,1)$.

This penalty function will establish the quality of a point $x$, considering the values of the functions $L_1$ and $L_2$. According to \cite{byrd1987trust}, the penalty parameter can decrease along the iterations if necessary but should be maintained whenever possible.

Finally, we present Algorithm \ref{alg:IRsimples}. This algorithm summarizes the ideas of the proposed method and is an almost direct adaptation of Algorithm 2 (``IR with globalization tests") in \cite{santos2022globalizaccao}.

\begin{algorithm} 
	\caption{Two-phase algorithm to minimize $L_1$ and $L_2$.}
	\label{alg:IRsimples}
	\begin{algorithmic}

\STATE{{\bf Step 0:} Initialization.
			
The parameters of the algorithm are defined as \( r \in (0,1) \), \( \theta_0 \in (0,1) \), and an initial estimate for the network weights \( (x^0) \). The iteration counter \( k \) is reset to zero.
}

\STATE{{\bf Step 1:} First Phase.
				
			If \( L_2(x^k) = 0 \), set \( y^k = x^k \). Otherwise, find a \( y^k \) that satisfies:
\begin{equation}
  \label{erres}
  L_2(y^k) \leq r L_2(x^k)
\end{equation}

\noindent

}
\STATE{	{\bf Step 2:} Penalty parameter.
		
If
\begin{equation} \label{penalidadetetakmais1}
  \Phi(y^{k},\theta_{k})\le \Phi(x^{k},\theta_{k})-\frac{1}{2} (1-r)^2 L_2(x^k),
\end{equation} 
}
do $\theta_{k+1}=\theta_k$.

Else, calculate:
\begin{equation} \label{eq:theta}
\theta_{k+1}=\frac{\frac{1}{2} (1-r^2)L_2(x^k) +rL_2(x^k) -L_2(y^k)    }{ L_1(y^k) - L_1(x^k) + L_2(x^k) -L_2(y^k)}.
\end{equation}
	
\STATE{	{\bf Step 3:} Second Phase.
	
Calculate \( d^k \in \mathbb{R}^n \) such that \( y^k + d^k \) satisfies:
\begin{equation} \label{testdirs}
L_1(y^k + d^k) \leq L_1(y^k) 
\end{equation} 
and
\begin{equation} \label{testmerits}
  \Phi(y^k + d^k, \theta_{k+1}) \leq  \Phi(x^k, \theta_{k+1}) -\frac{1}{2}(1-r)^2 L_2(x^k).  
\end{equation}

}
	
\STATE{	{\bf Step 4:} Actualization
	
	Set \( x^{k+1} := y^k + d^k \). If a stopping criterion is not satisfied, increment \( k \) and return to Step 1.
}
		
	\end{algorithmic}
\end{algorithm}

In Step 2, a specific algorithm introduced and discussed in \cite{santos2022globalizaccao}, is used to determine \( \theta_k \). By defining \( \theta_k \) in this way, we can ensure that \( \{\theta_k\}_k \) will be a non-increasing sequence in \( (0,1) \), which is crucial for the discussion of the algorithm's convergence. This also allows us to show that there exists some \( d^k \) satisfying the conditions \eqref{testdirs} and \eqref{testmerits}, i.e., that the second phase of the algorithm is well defined.

Unlike its predecessor, Algorithm \ref{alg:IRsimples} presented here does not consider the step size \( \|d^k\| \) in the descent condition \eqref{testdirs}. We made this choice because implementing this condition would imply a higher computational effort at each algorithm step. This adaptation is equivalent to fixing the parameter \( \sigma = 0 \) in the original algorithm, making the descent condition \eqref{testdirs} presented here weaker than its predecessor. Therefore, the theoretical convergence guarantee results presented in \cite{santos2022globalizaccao} are not maintained. However, more general versions of this algorithm have already been explored (see \cite{moreno2023globalizaccao} for the study of RI methods with non-monotonic descent conditions) with theoretical guarantees, which makes this version of the algorithm, even without a well-established theoretical guarantee, a viable option to be explored.

Note that the proposed method is a general-purpose algorithm since it does not specify which method should be used to estimate the new iterate satisfying the descent conditions in each phase. The ideal strategy, as happens in IR, is an entire research area.

Another critical factor in the performance of Algorithm \ref{alg:IRsimples} is the choice of parameters \( \theta_0 \) and \( r \). As in IR, the parameter \( r \) determines the requirement level of the first phase. The smaller \( r \) is, the greater the expected decrease of the function \( L_2 \) value already in the first phase, and the descent condition \eqref{testmerits} becomes more demanding. The ideal choice of these parameters is linked to the problem to be solved; however, according to the numerical tests carried out in \cite{moreno2023globalizaccao} and \cite{santos2022globalizaccao}, it is expected that values of \( r \) and \(  \theta_0 \) close to 1 will generate good results.

\section{Experiments and Results}

As discussed earlier, the optimization algorithm proposed in this work was developed for training networks whose performance depends on two objectives represented by different terms in the loss function. And, as outlined in \cite{raissi2019physics}, this is the case sometimes when using a PINN to estimate the solution of a differential equation.

For this reason, we will consider testing the performance of Algorithm \ref{alg:IRsimples} on Physics Informed Neural Networks. In particular, in this work, we will test the performance of Algorithm \ref{alg:IRsimples} in training PINN to solve two different problems:

\subsection{PINN to solve Burgers' equation}

First we will train a PINN aimed at solving the Burgers' differential equation with viscosity, sinusoidal initial condition, and Dirichlet boundary conditions, namely:

\begin{eqnarray*}
    &&u_t + u u_x - (0.01/\pi)u_{xx} = 0, \quad x \in [-1,1], \quad t \in [0,1], \\
    &&u(0,x) = -\sin(\pi x), \\
    &&u(t,-1) = u(t,1) = 0.
\end{eqnarray*}

This problem has become a standard benchmark in testing Physics Informed Neural Networks because the solution exhibits a discontinuity, making it a critical case for study.

To train a Neural Network that approximates the solution to the above problem, we have two loss terms involved: one is responsible for ensuring that the initial and boundary conditions are satisfied, and the other is related to satisfying the differential equation. The function responsible for satisfying the boundary and initial conditions is traditionally an MSE evaluated over a set of fixed points on the domain boundary, called data points. Following the notation used in \cite{raissi2019physics}, we denote this function by $MSE_u$. Similarly, to assess the satisfaction of the differential equation, collocation points are generated in the domain, and we estimate the differential equation at those points. The loss function associated with this objective is also traditionally an MSE, denoted by $MSE_f$. Thus, in this case, the aim is to minimize these two functions together. 

In the following section, we will compare the performance of a Neural Network trained using the simple sum loss function $MSE_u + MSE_f$, as proposed in \cite{raissi2019physics}, with another network with the same architecture trained using Algorithm \ref{alg:IRsimples}.

In all cases, we used a fully connected deep Neural Network with nine layers of 20 neurons each, with a hyperbolic tangent activation function and an output layer with one neuron without activation function. We randomly generated 100 training points on the boundary to ensure that the initial and boundary conditions are satisfied and 10.000 collocation points using Latin Hypercube Sampling (LHS) as in \cite{raissi2019physics}.

\subsubsection{Reference example}

We consider as a reference example the network trained using the code available at \url{https://github.com/314arhaam/burger-pinn}. In this case, training is carried out using the Adam optimizer. We consider this because it is one of the most commonly used optimizers in PINN training to date. This will serve as a baseline example for making comparisons with the implementation of Algorithm \ref{alg:IRsimples}, which was also developed using this code as a base with some modifications.

\begin{figure}[!t]
\centering
\subfloat[Loss functions.]{\includegraphics[width=1.45in]{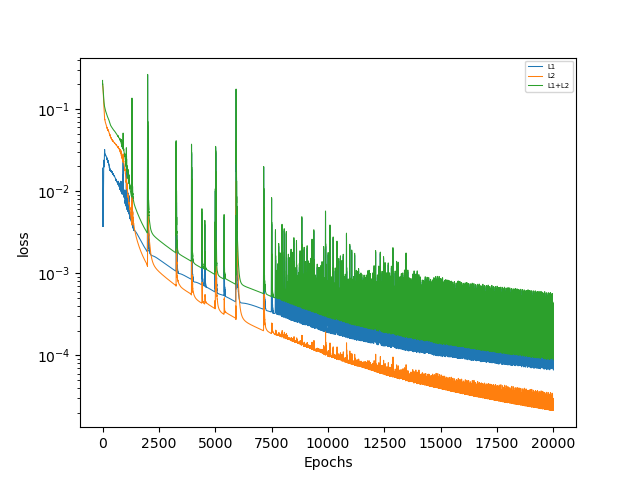}%
\label{original_custo_7}}
\subfloat[Estimated solution $u(t,x)$.]{\includegraphics[width=1.9in]{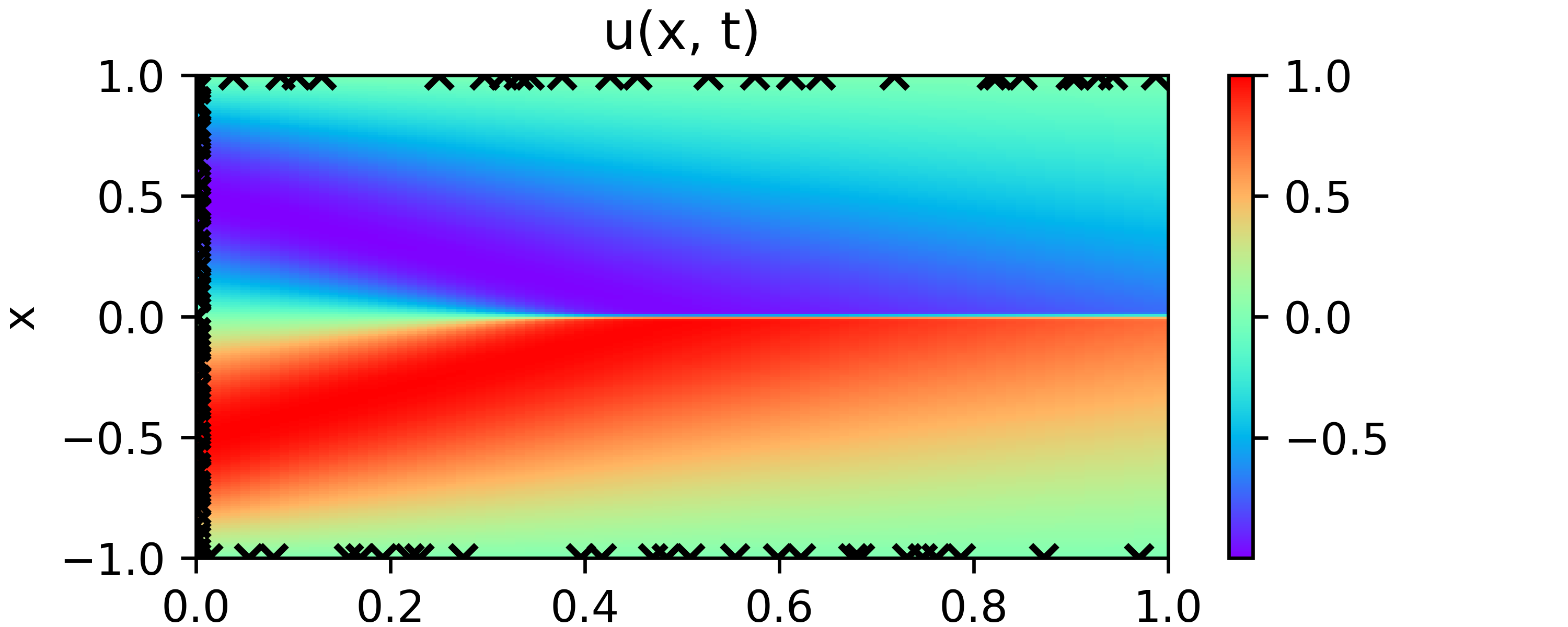}%
\label{original_campo_7}}

\subfloat[Estimated solution $u(t,x)$ vs real solution.]{\includegraphics[width=2.9in]{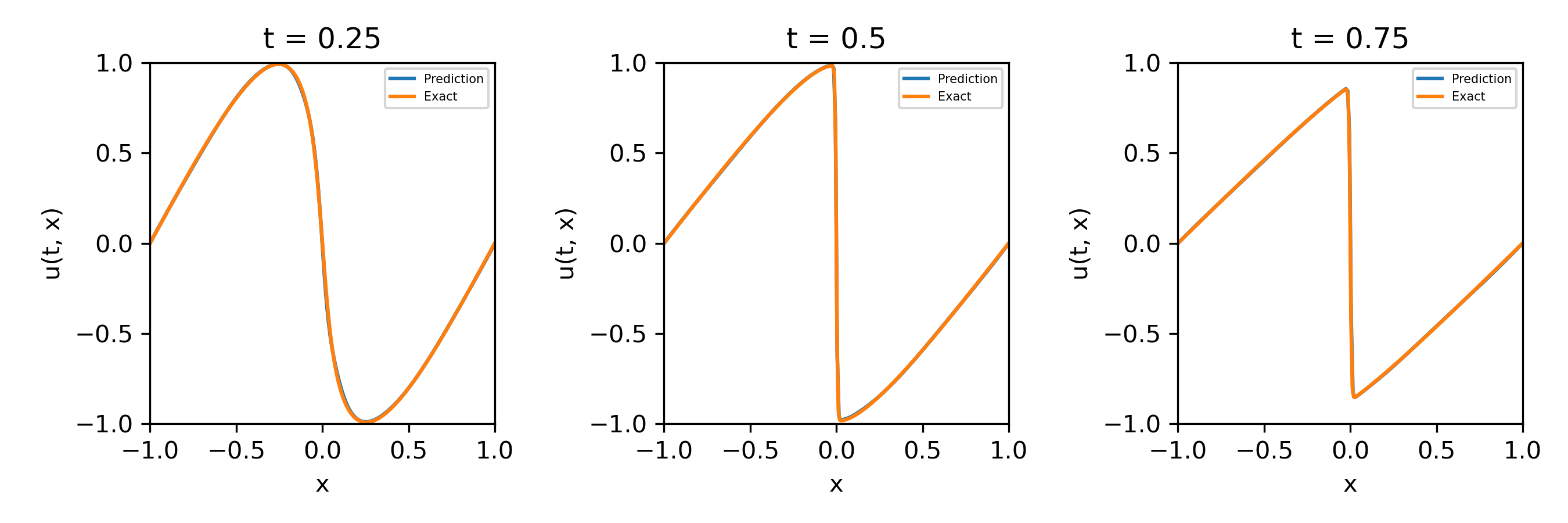}%
\label{original_curvas_7}}
\caption{Results for the Neural Network trained to solve Burgers' equation with the Adam optimizer. Left: History of the functions $L_1=MSE_f$, $L_2=MSE_u$ and $L_1+L_2$. Right: Solution $u(t,x)$ estimated by the Neural Network. Bottom: Solution obtained by the Neural Network vs real solution of the Burgers' equation for specific times. }
\label{fig_original_custo_campo_curvas_7}
\end{figure}

In the original code, network training is done using a constant learning rate of $0.0005$, with the loss function $MSE_u + MSE_f$. In Figure \ref{original_custo_7} and \ref{original_campo_7}, we show the behavior of $L_1, L_2$ and $L_1+L_2$ during training with 20.000 epochs and the results obtained by the network. In Figure \ref{original_curvas_7}, we show the predicted solutions at specific time instances $t=0.25$, $t=0.5$, and $t=0.75$ in contrast to the real solution, which is available in \url{https://github.com/maziarraissi/PINNs}.

\begin{figure}[!t]
\centering
\subfloat[3000 epochs.]{\includegraphics[width=2.9in]{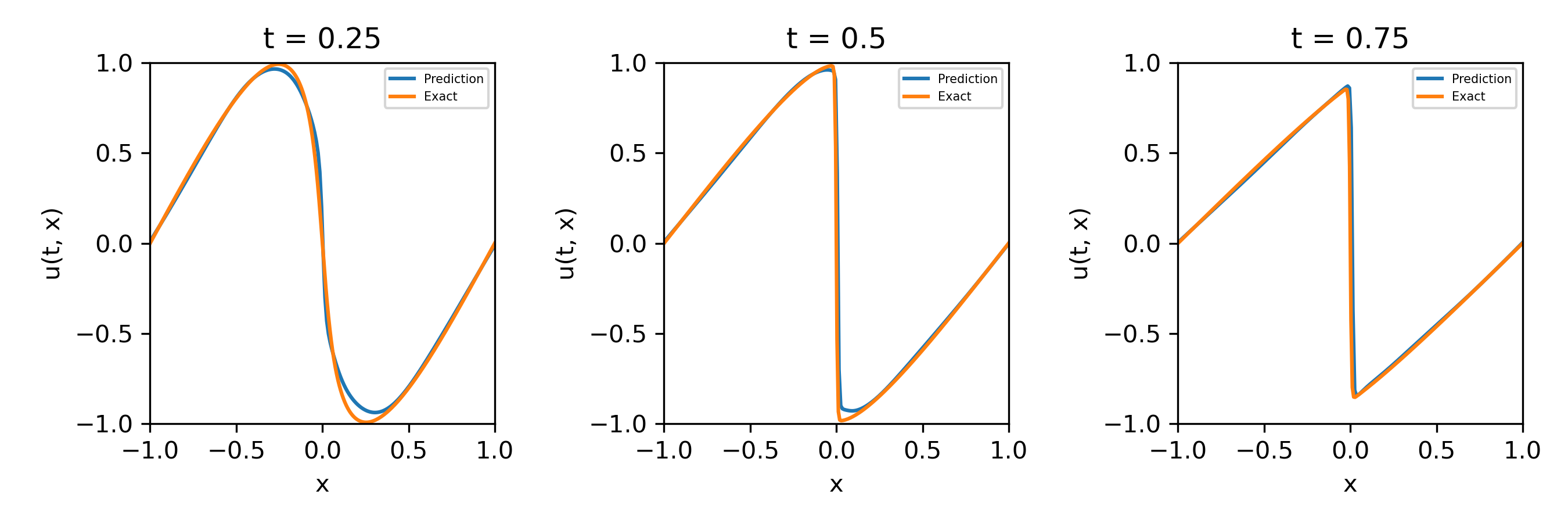}%
\label{7_original_curvas_t1}}

\subfloat[5000 epochs.]{\includegraphics[width=2.9in]{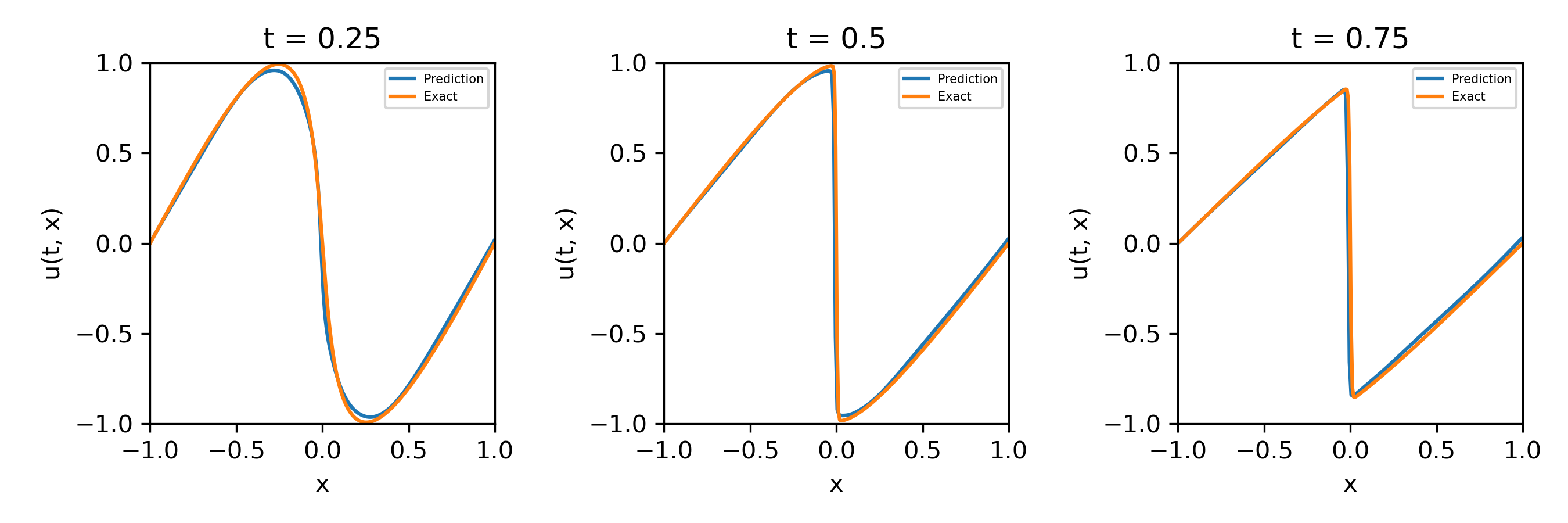}%
\label{7_original_curvas_t2}}
\caption{Solution obtained by the Neural Network vs real solution of the Burgers' equation for a Neural Network trained with the Adam optimizer with 3000 and 5000 epochs. }
\label{fig_original_etapas_7_curvas}
\end{figure}

To better understand the convergence speed, in Figures \ref{7_original_curvas_t1} and \ref{7_original_curvas_t2}, we display the predicted solutions at three specific times for the same network trained with 3000 and 5000 epochs, respectively.

\begin{figure}[!t]
\centering
\subfloat[1980 epochs.]{\includegraphics[width=2.9in]{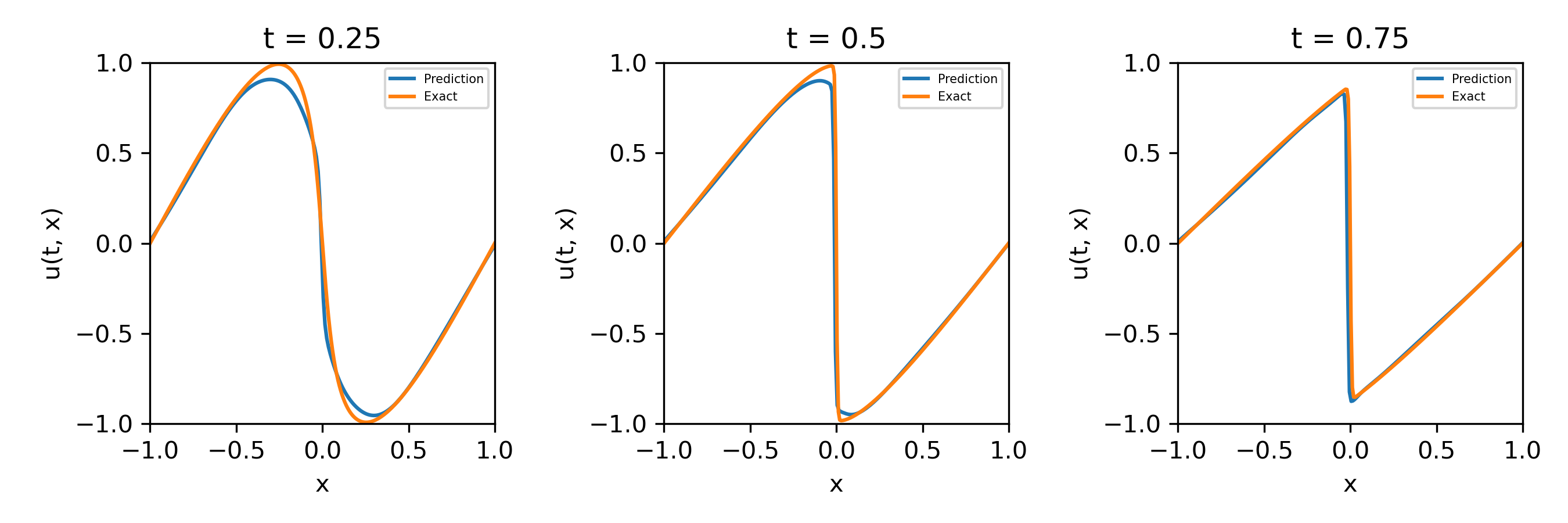}%
\label{comparacao_1}}

\subfloat[2000 epochs.]{\includegraphics[width=2.9in]{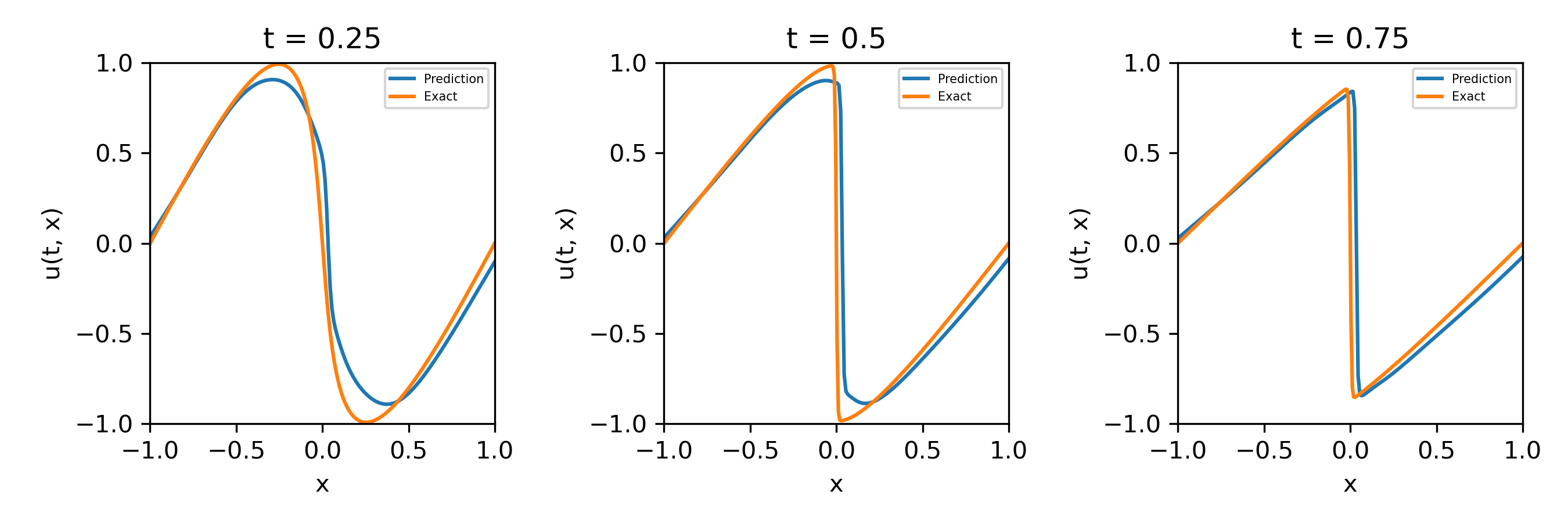}}%
\label{comparacao_2}
\caption{Solution obtained by the Neural Network vs real solution of Burgers' equation for Neural Network trained with the Adam optimizer after 1980 and 2000 epochs.}
\label{fig_orig_imag_2_curvas_comp}
\end{figure}

Sudden oscillations are observed in the loss functions during the training process. This means that increasing the number of epochs does not necessarily result in better loss function values or improved network performance. For example, the approximation made by the network trained after the 1980 epochs shows better results than the same network after the 2000 epochs during the same training process (Figure \ref{fig_orig_imag_2_curvas_comp}).

In Algorithm \ref{alg:IRsimples}, we expect to control this phenomenon by imposing descent conditions.

Although the oscillatory behavior of the loss function is undesirable, the method ends up converging reasonably well. This is evident in the results presented in Figure \ref{fig_original_custo_campo_curvas_7} for the network trained after 20.000 epochs.

\subsubsection{The Proposed Method}

We used the reference code as a base for implementing Algorithm \ref{alg:IRsimples}, with the necessary modifications.

One of the main characteristics of Algorithm \ref{alg:IRsimples} is the imposition of descent conditions in each phase. The objective in each phase is to find a new iterate that satisfies these conditions. This work does not extensively discuss the method used to find these new iterates. However, while the descent conditions are designed to be theoretically feasible, meaning that it should always be possible to find an iterate that satisfies them, in practice, many attempts may be required. Therefore, it is evident that each complete iteration of Algorithm \ref{alg:IRsimples} may require more time and computational effort than an iteration of a simpler optimization algorithm, such as Adam.

For this reason, for a fair comparison with other search algorithms, we count the number of partial iterations and the number of epochs. Each search attempt in a direction is considered a partial iteration. On the other hand, we consider that an epoch has occurred once all phases of the algorithm have been completed. Specifically, the number of epochs refers to the counter $k$ in Algorithm \ref{alg:IRsimples}.

 To limit the computational effort in each phase, we introduce a crucial element, the hyperparameter $ it_{ {max}} $. This parameter, which will be fixed in each case, limits the number of internal iterations within each phase of the algorithm.

In the following tests, we consider the Adam optimizer with the same learning rate as the original code as a helper to determine the descent directions within each phase.

The choice of parameters $\theta_0=0.8$ and $r=0.99$ was made following the notions discussed in \cite{moreno2023globalizaccao}.

\begin{itemize}
    \item First set of tests: For this set of tests, we will train the Neural Network using Algorithm \ref{alg:IRsimples}, fixing $L_1 = \text{MSE}_f$ and $L_2 = \text{MSE}_u$.

The number of internal iterations in each phase is limited to $it_{\text{max}} = 150$. To minimize the function $L_2$ in the first phase, we apply the Adam optimizer using the gradient information of $L_2$. However, in the second phase, we aim to minimize not only $L_1$ but also the penalty function $\Phi$. Therefore, we apply the optimizer with the gradient information of the linear combination $\alpha L_1 + \beta L_2$, where $\alpha$ and $\beta$ are hyperparameters that have been determined empirically in this work. For this set of tests, we fix them at $\alpha = 0.5$ and $\beta = 4$. 

In Figure \ref{fig_Testes_1_1_custo_campo_curvas}, we show the behavior of the loss functions ($L_1,L_2$ and $L_1+L_2$) during training with 20.000 internal iterations, resulting in 982 epochs, along with the estimated solution by the Neural Network for Burgers' equation in the established domain after complete training, and at three particular time instances ($t$). Notably, the behavior of the loss function is more stable, as expected from the imposed descent conditions. Furthermore, after 20.000 partial iterations, the solution estimated by the network overlaps (at least at the studied time instances $t$) with the real solution.

\begin{figure}[!t]
\centering
\subfloat[Loss functions.]{\includegraphics[width=1.45in]{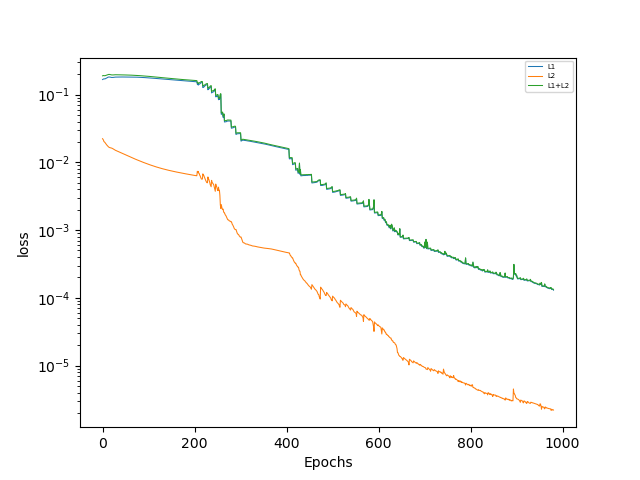}%
\label{Testes_1_1_custo}}
\subfloat[Estimated solution $u(t,x)$.]{\includegraphics[width=1.9in]{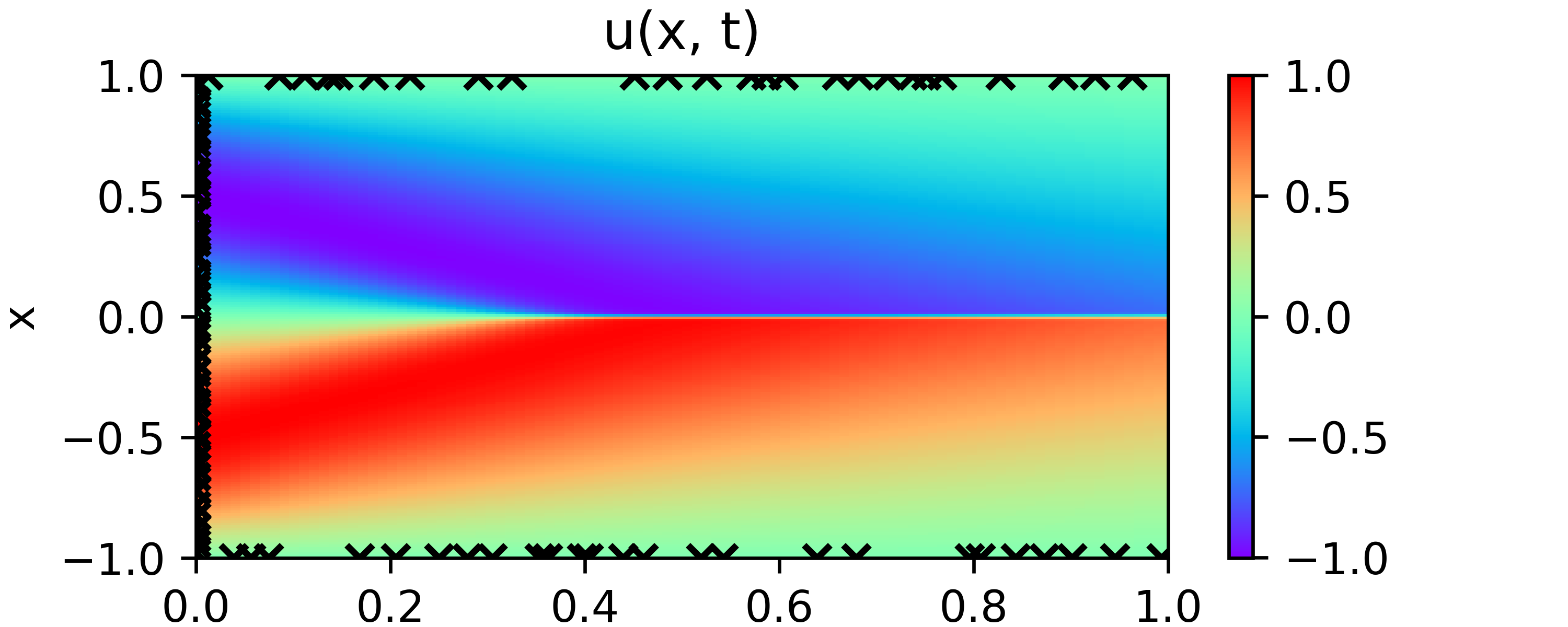}%
\label{Testes_1_1_campo}}

\subfloat[Estimated solution $u(t,x)$ vs real solution.]{\includegraphics[width=2.9in]{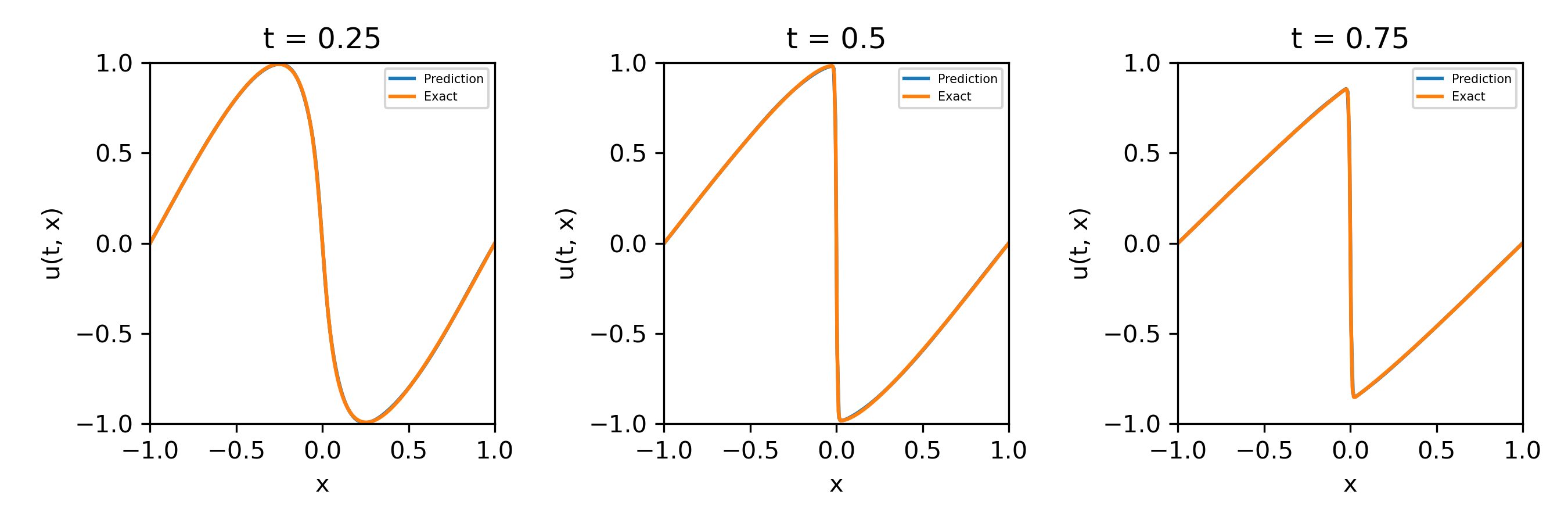}%
\label{Testes_1_1_curvas}}
\caption{Results for the Neural Network trained following the Algorithm \ref{alg:IRsimples} considering $L_1=MSE_f$ and $L_2=MSE_u$. Left: History of functions $L_1$, $L_2$ and $L_1+L_2$. Right: Solution $u(t,x)$ estimated by the Neural Network. Bottom: Solution obtained by the Neural Network vs real solution of Burgers' equation for specific moments of time. }
\label{fig_Testes_1_1_custo_campo_curvas}
\end{figure}

To better understand the convergence behavior, in Figure \ref{Testes_1_1}, we show the comparison of the real solution of Burgers' equation for different time values with the solution estimated by the Neural Network at different stages of training.

\begin{figure}[!t]
\centering
\subfloat[3000 Internal iterations.]{\includegraphics[width=2.9in]{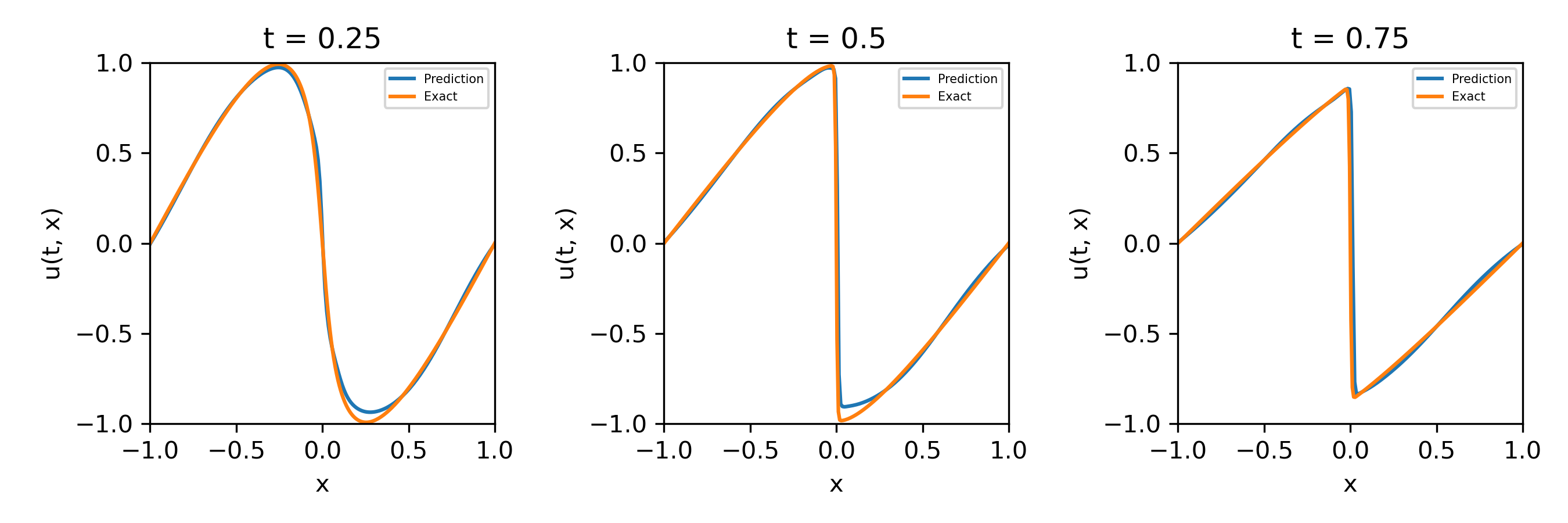}%
\label{Testes_1_1_t1}}

\subfloat[5000 Internal iterations.]{\includegraphics[width=2.9in]{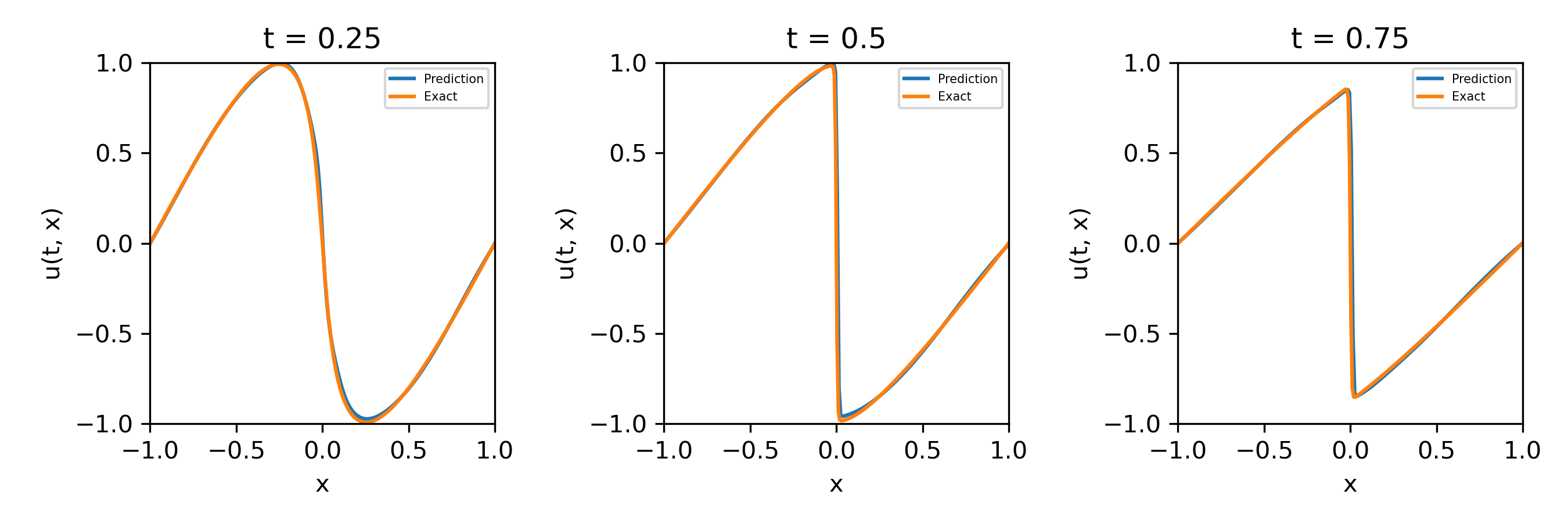}%
\label{Testes_1_1_t2}}
\caption{Solution obtained by the Neural Network vs real solution of Burgers' equation for Neural Networks trained with the Algorithm \ref{alg:IRsimples} where $L_1=MSE_f$ and $L_2=MSE_u$ with 3000 and 5000 internal iterations.}
\label{Testes_1_1}
\end{figure}

Comparing Figures \ref{fig_original_etapas_7_curvas} and \ref{Testes_1_1}, which show the evolution of the network's estimates at two different epochs in contrast to the real solution, we notice that in the studied epochs, the proposed method showed a faster convergence than the one trained with Adam. However, as is evident when studying the results of the fully trained networks (Figures \ref{original_curvas_7} and \ref{Testes_1_1_curvas}), both methods converge.

\item Second set of tests: To show that the convergence of Algorithm \ref{alg:IRsimples} is independent of the choice of which term of the loss function is associated with which function, we consider, contrary to the previous test set, $L_1 = \text{MSE}_u$ and $L_2 = \text{MSE}_f$.

We set the maximum number of internal iterations to $it_{\text{max}} = 250$. And, as in the previous case, in the second phase, we use the Adam optimizer with gradient information of $\alpha L_1 + \beta L_2$, where now $\alpha = 2$ and $\beta = 1.5$.

In Figure \ref{Testes_1_2_custo}, we show the behavior of the loss functions after 20.000 internal iterations, which represents only 492 epochs, or in other words, 492 complete iterations of Algorithm \ref{alg:IRsimples}. In Figure \ref{Testes_1_2_campo}, we show Burgers' equation solution prediction made by the same network across the domain. In contrast, in Figure \ref{Testes_1_2_curvas}, this same solution is observed compared to the real solution at three specific times.

Notice that the values of the loss functions may show some oscillations but not as intense as in the reference example. This occurs because within the phases of Algorithm \ref{alg:IRsimples}, an adequate method is not used. Adequate methods use region of trust strategies, as is commonly used in IR algorithms (\cite{bellavia2020inexact}, \cite{bellavia2023stochastic}, \cite{birgin2005local}, \cite{martinez2005inexact}), which allow ensuring a non-increasing behavior of the function to be minimized throughout the optimization process. Thus, the maximum number of iterations within the phases ($it_{\text{max}}$) is often reached, allowing for this oscillatory behavior of the loss functions. However, the method shows convergence at the end of training, as the estimate made by the network overlaps with the real solution.

\begin{figure}[!t]
\centering
\subfloat[Loss functions.]{\includegraphics[width=1.45in]{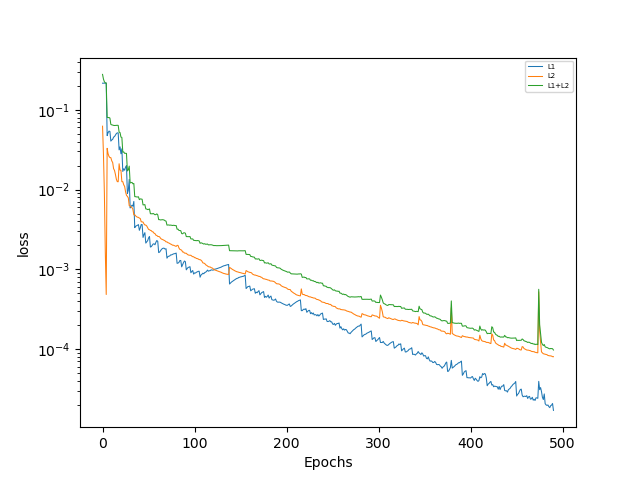}%
\label{Testes_1_2_custo}}
\subfloat[Estimated solution $u(t,x)$.]{\includegraphics[width=1.9in]{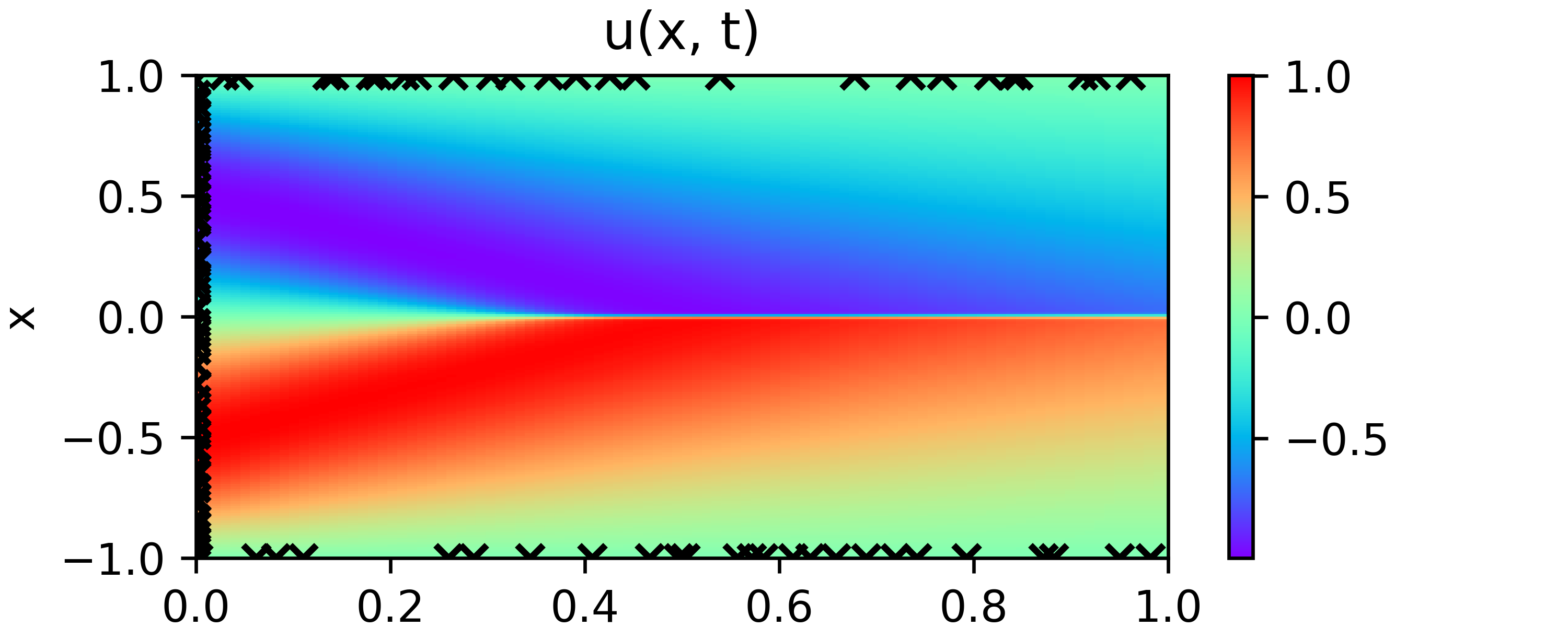}%
\label{Testes_1_2_campo}}

\subfloat[Estimated solution $u(t,x)$ vs real solution.]{\includegraphics[width=2.9in]{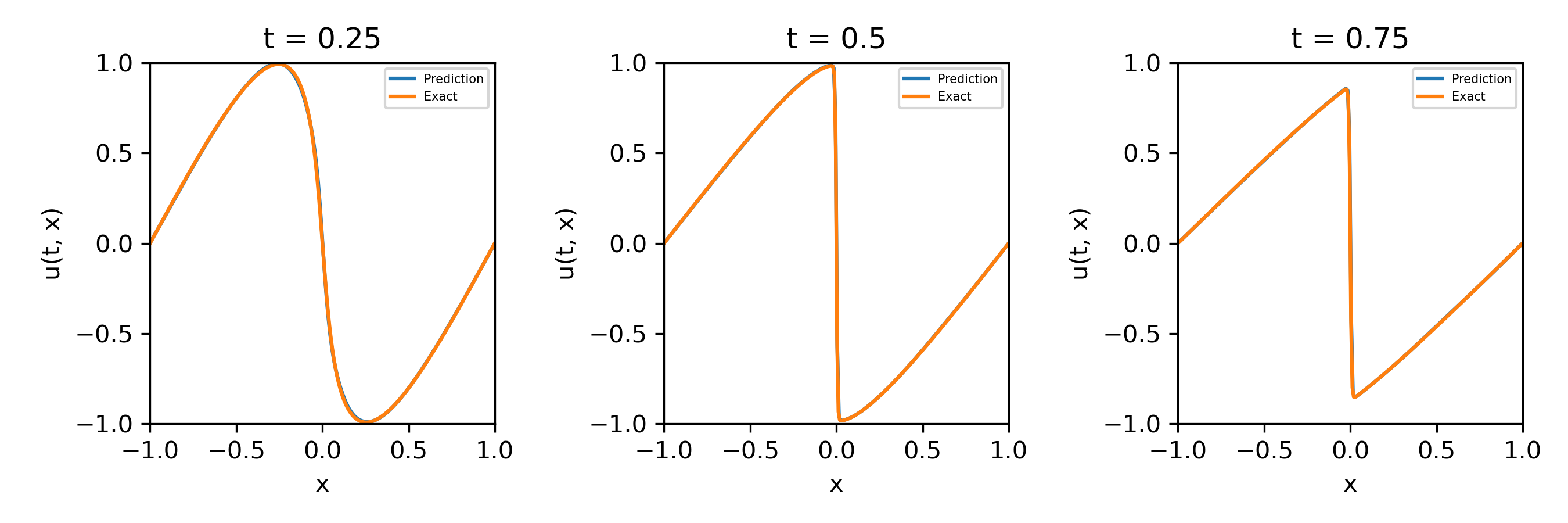}%
\label{Testes_1_2_curvas}}
\caption{Results for the Neural Network trained following the Algorithm \ref{alg:IRsimples} considering $L_1=MSE_u$ and $L_2=MSE_f$. Left: History of functions $L_1$, $L_2$ and $L_1+L_2$. Right: Solution $u(t,x)$ estimated by the Neural Network. Bottom: Solution obtained by the Neural Network vs real solution of the Burgers' equation for specific times.}
\label{fig_Testes_1_2_custo_campo_curvas}
\end{figure}

In Figure \ref{Testes_1_2}, we show the network's estimate of Burgers' equation solution at three times for the network trained with 3000 and 5000 partial iterations.

\begin{figure}[!t]
\centering
\subfloat[3000 Internal iterations.]{\includegraphics[width=2.9in]{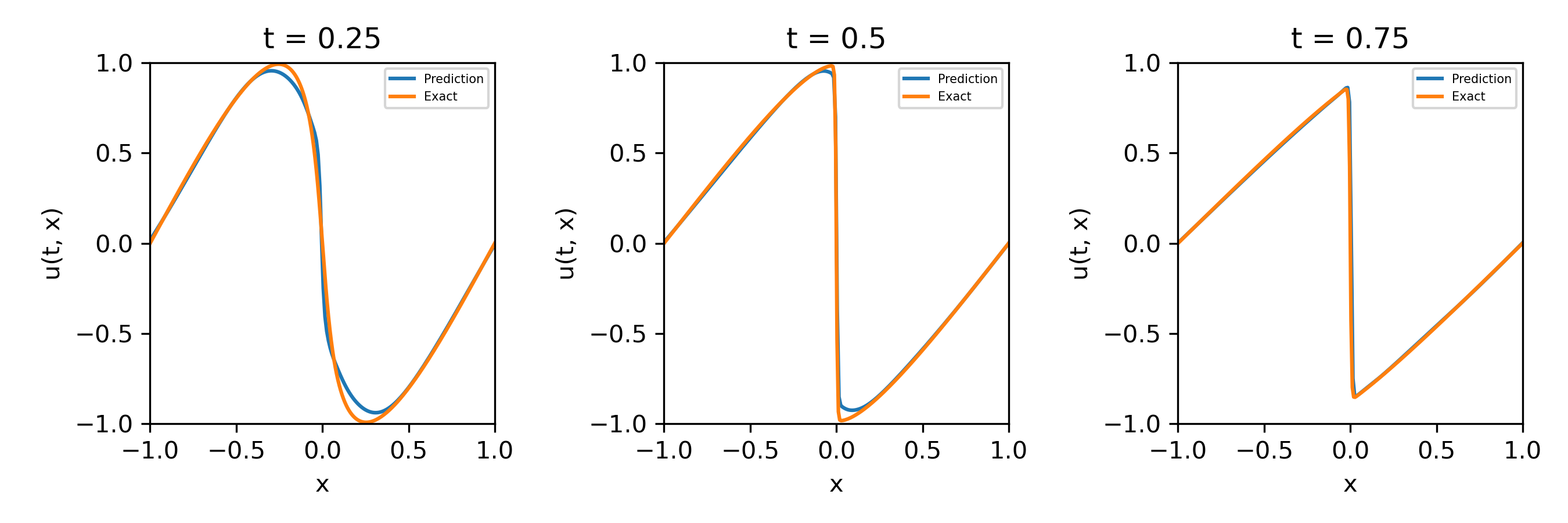}%
\label{Testes_1_2_t1}}

\subfloat[5000 Internal iterations.]{\includegraphics[width=2.9in]{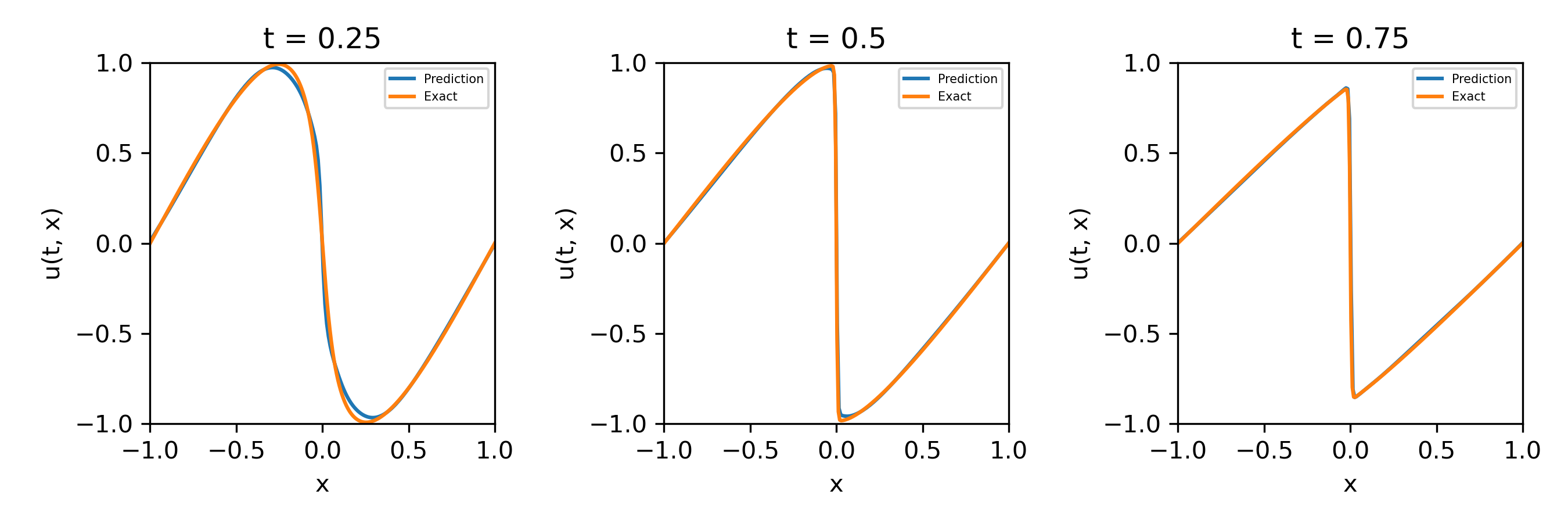}%
\label{Testes_1_2_t2}}
\caption{Solution obtained by the Neural Network vs real solution of  Burgers' equation for  a Neural Network trained with the Algorithm \ref{alg:IRsimples} where $L_1=MSE_u$ and $L_2=MSE_f$ with 3000 and 5000 internal iterations.}
\label{Testes_1_2}
\end{figure}

Comparing the evolution of the results obtained by Algorithm \ref{alg:IRsimples} at two moments of training (Figures \ref{fig_Testes_1_2_custo_campo_curvas} and \ref{Testes_1_2}) with the results obtained by the reference example at these same moments (Figures \ref{fig_original_custo_campo_curvas_7} and \ref{fig_original_etapas_7_curvas}), we conclude that the convergence of the method proposed here was, at least, comparable to that of the reference algorithm.

\end{itemize}

We will develop another example in the next section to demonstrate the competence of the proposed method in various scenarios.

\subsection{PINN to solve Heat equation}

Here, our objective will be to solve the one-dimensional heat equation with a constant thermal coefficient and without a heat source, namely:

\begin{eqnarray*}
    &&u_t= k u_{xx}, \quad x \in [0,L], \quad t \in [0,5], \\
    &&u(0,x) = \sin(\pi x/L)+ \sin(3 \pi x/L), \\
    &&u(t,0) = u(t,L) = 0.
\end{eqnarray*}

To solve this problem, we will use the parameters and network structure suggested in \cite{popa2022solving}. Thus, we set $k=1$ and $L=10$ in the equations above. We also consider an MLP with four hidden layers and 50 neurons each, and $tanh$ as the activation function. The collocation points are generated as in the previous example.

Following a strategy similar to the previous section, we will first present the results obtained for the network trained using traditional methods and then compare them with those obtained using the Algorithm \ref{alg:IRsimples}.

\subsubsection{Reference example} 

We trained the network using the Adam optimizer with a learning rate of $0.0005$.

Again, we consider the loss function to be the sum of MSE functions associated with the satisfaction of the PDE ($MSE_f$) and the satisfaction of the initial and boundary conditions ($MSE_u$).

\begin{figure}[!t]
\centering
\subfloat[Loss functions.]{\includegraphics[width=1.45in]{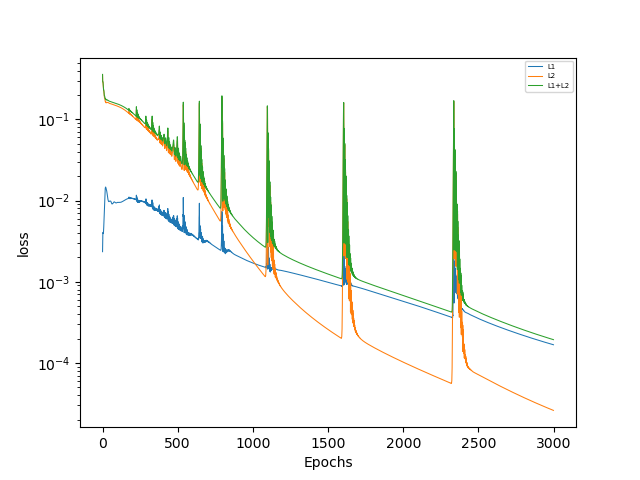}%
\label{original_custo_16}}
\subfloat[Estimated solution $u(t,x)$.]{\includegraphics[width=1.9in]{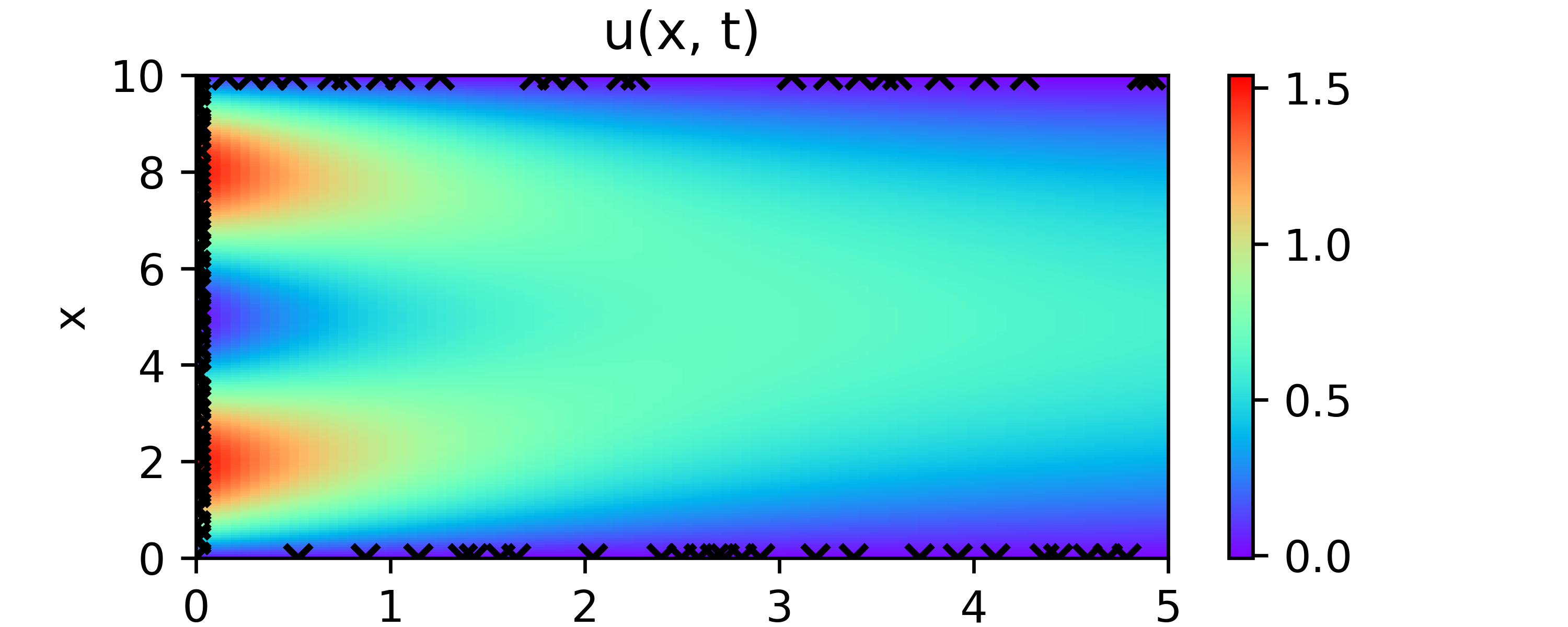}%
\label{original_campo_16}}

\subfloat[Estimated solution $u(t,x)$ vs real solution.]{\includegraphics[width=2.9in]{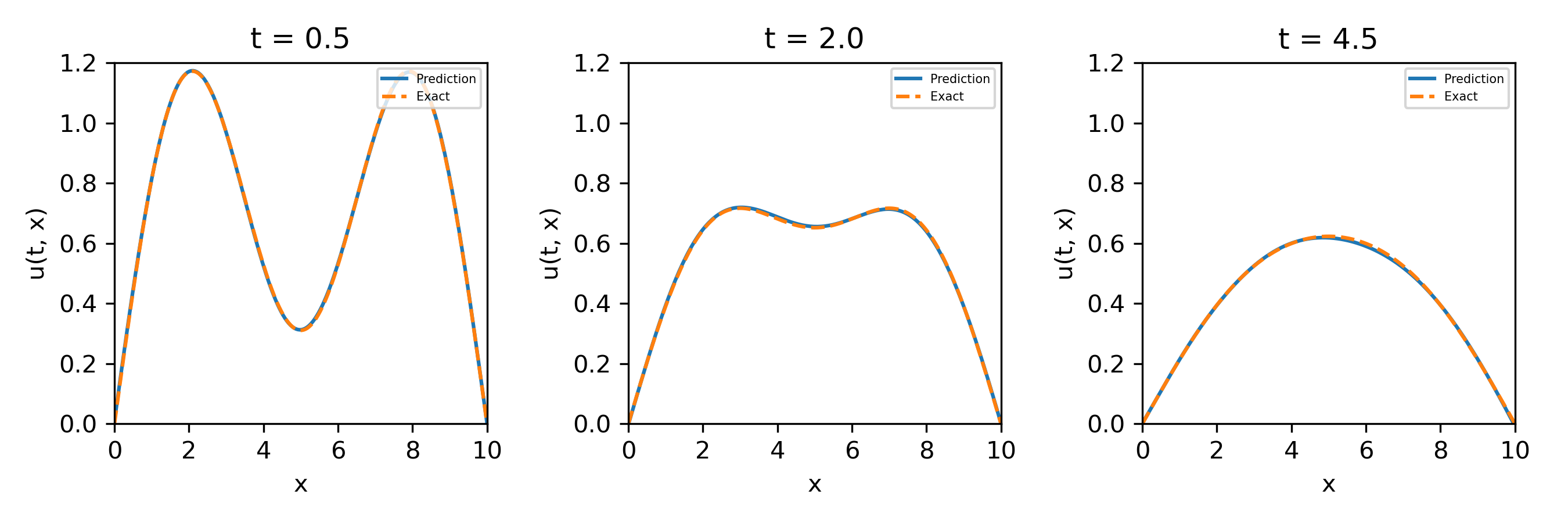}%
\label{original_curvas_16}}
\caption{Results for the Neural Network trained with the Adam optimizer. Left: History of the functions $L_1=MSE_f$, $L_2=MSE_u$ and $L_1+L_2$. Right: Solution $u(t,x)$ estimated by the Neural Network. Bottom: Solution obtained by the Neural Network vs real solution of the Heat equation for specific times. }
\label{fig_original_custo_campo_curvas_16}
\end{figure}

In Figure \ref{fig_original_custo_campo_curvas_16}, we show the results obtained after 3000 epochs.

\begin{figure}[!t]
\centering
\subfloat[2320 epochs.]{\includegraphics[width=2.9in]{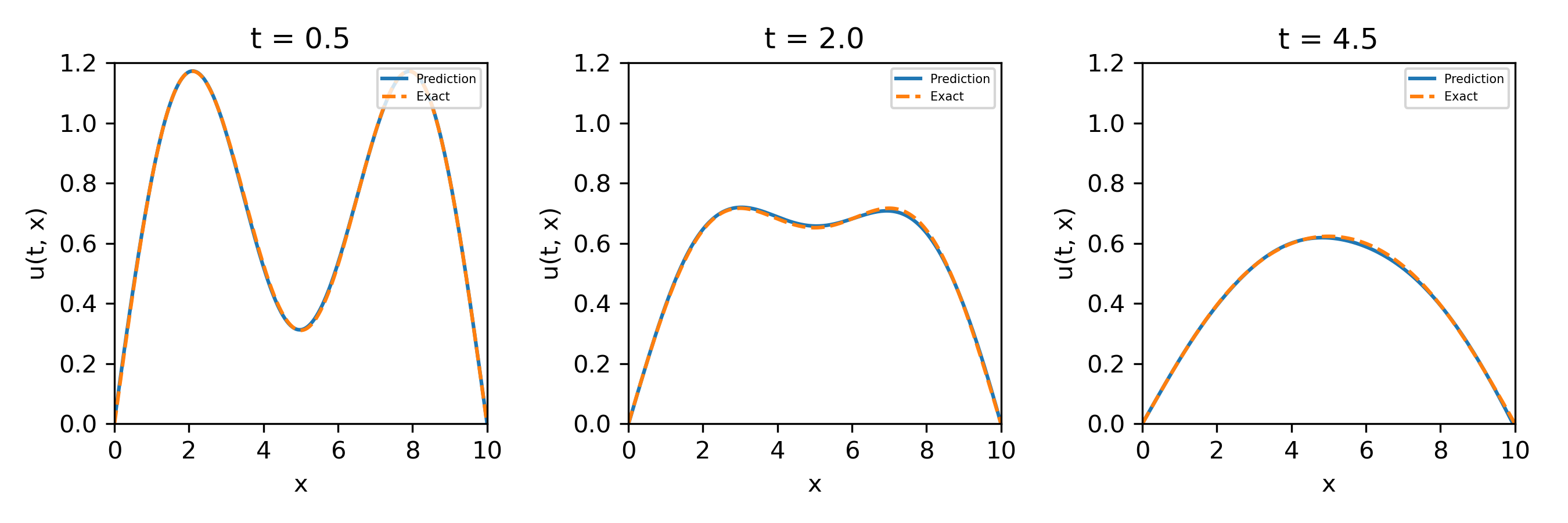}%
\label{comparacao_2_1}}

\subfloat[2340 epochs.]{\includegraphics[width=2.9in]{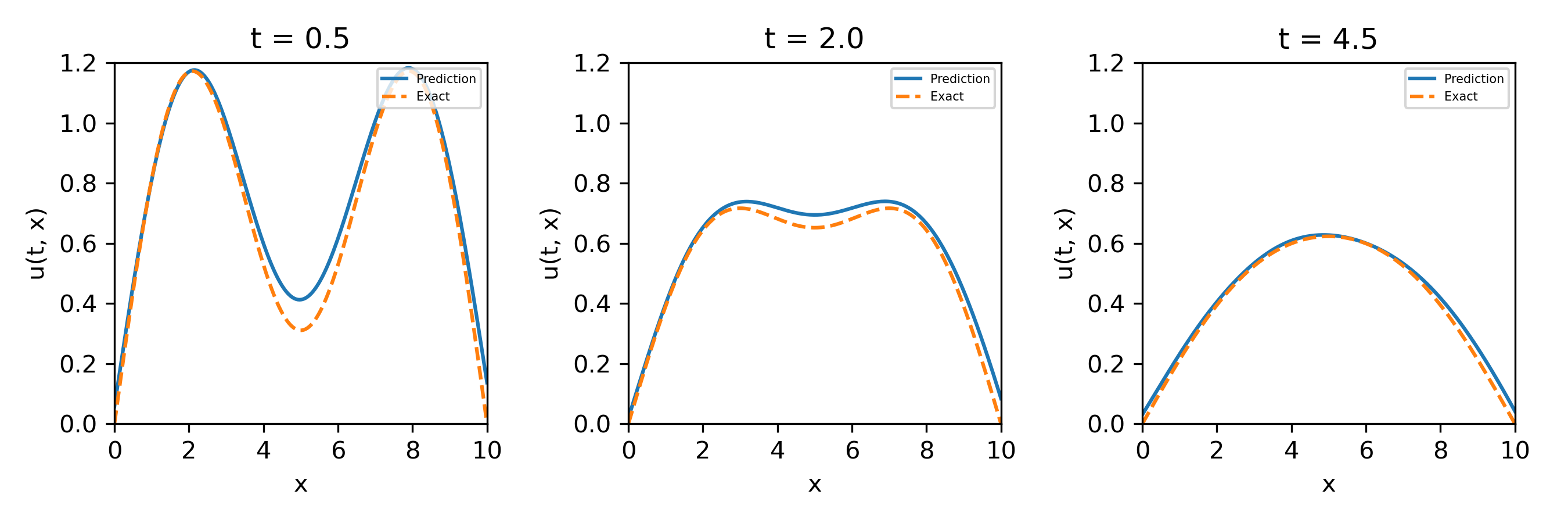}%
\label{comparacao_2_2}}
\caption{Solution obtained by the Neural Network vs real solution of Heat equation for Neural Network trained with the Adam optimizer after 2320 and 2340 epochs.}
\label{fig_orig_imag_16_curvas_comp}
\end{figure}

In this experiment, the phenomenon observed in the previous example repeats: the loss function shows strong oscillations. This phenomenon causes the results to vary depending on the epoch at which the training is stopped, potentially leading to worse results. This is evident when we observe, for example, the network's estimated result after  2320 and 2340 epochs (Figure \ref{fig_orig_imag_16_curvas_comp}). In this case, after achieving a reasonable result in epoch 2320, 20 epochs later the network's estimation significantly deteriorates in quality.

The following section presents the results obtained with Algorithm \ref{alg:IRsimples}.

\subsubsection{The Proposed Method}

We will use a strategy similar to the one discussed for the proposed method in the previous experiments. We will use the Adam optimizer as an auxiliary and limit the number of internal iterations in each phase. 

\begin{itemize}
    \item First set of tests: The maximum number of iterations within each phase is $it_{max}=100$. On the other hand, the Adam optimizer within phase 2 uses information from the gradient of $\alpha L_1+ \beta L_2$ where $\alpha=4$ and $\beta=1$.
    
    In Figure  \ref{fig_Testes_2_1_custo_campo_curvas} we show the loss function behaviour and the estimated solution.

    \begin{figure}[!t]
        \centering
        \subfloat[Loss functions.]{\includegraphics[width=1.45in]{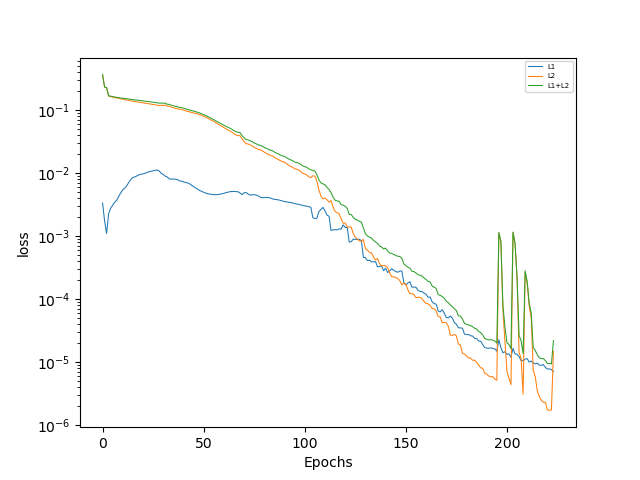}%
        \label{Testes_2_1_custo}}
        \subfloat[Estimated solution $u(t,x)$.]{\includegraphics[width=1.9in]{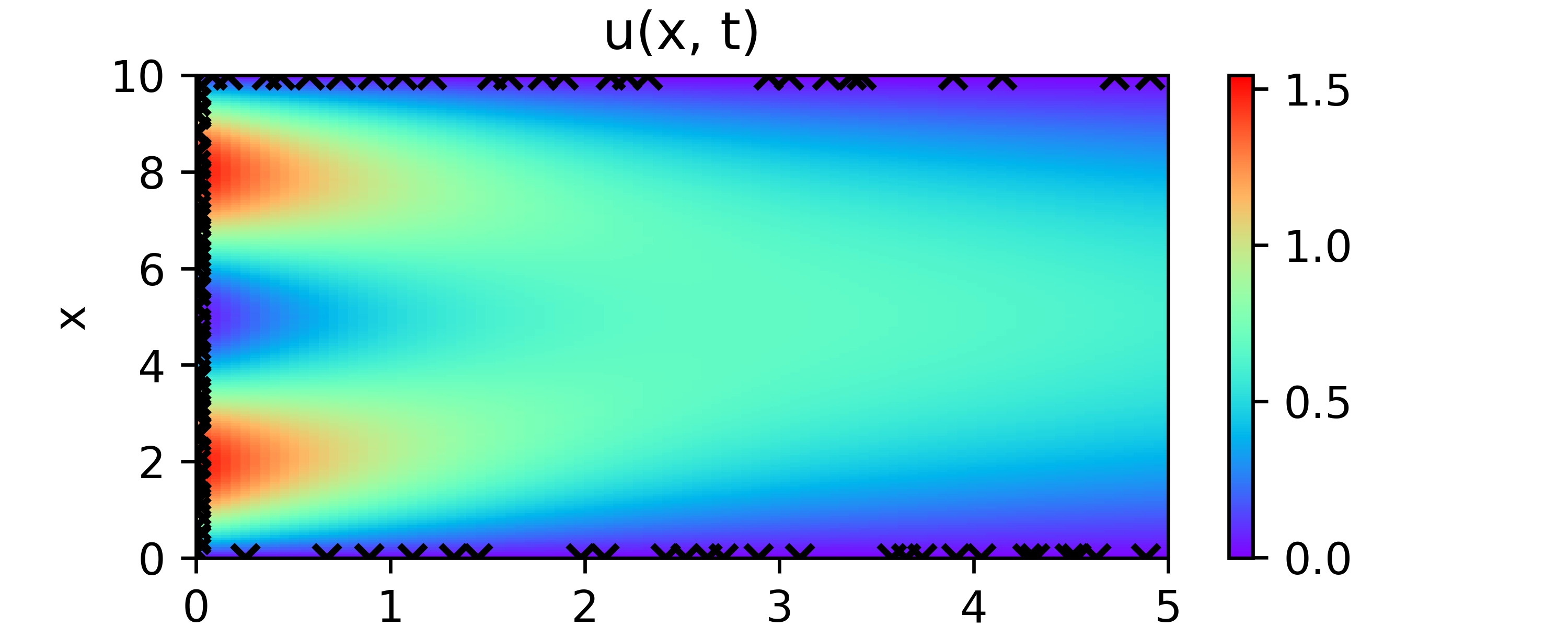}%
        \label{Testes_2_1_campo}}
        
        \subfloat[Estimated solution $u(t,x)$ vs real solution.]{\includegraphics[width=2.9in]{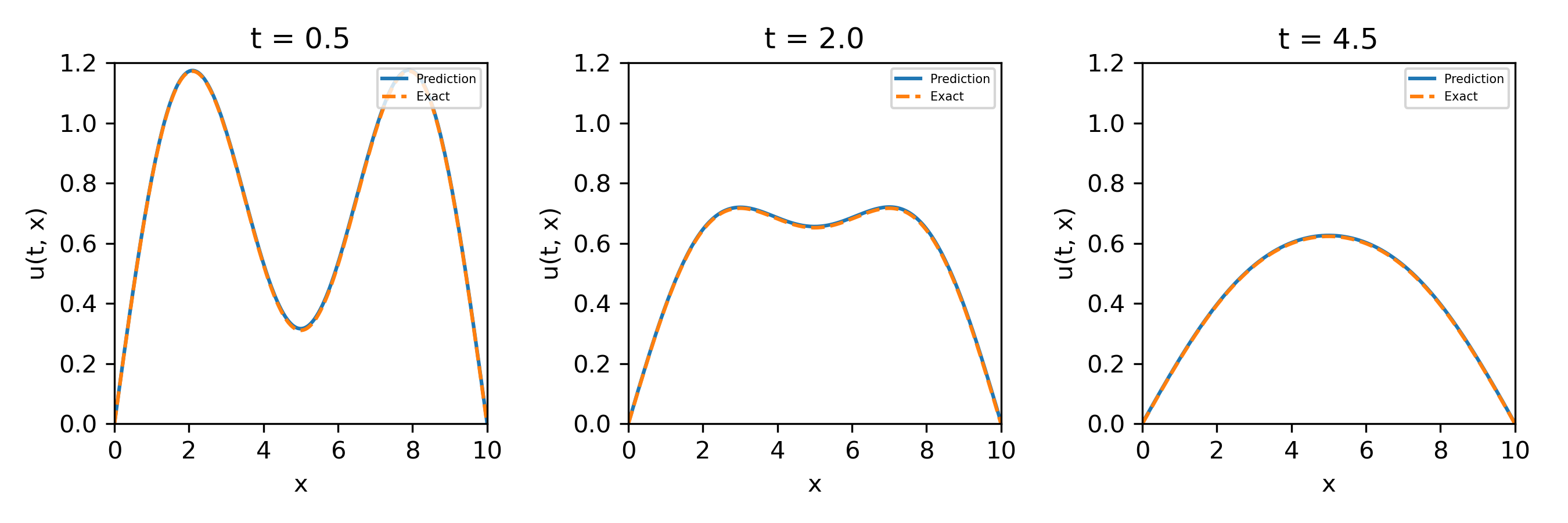}%
        \label{Testes_2_1_curvas}}
        \caption{Results for the Neural Network trained following the Algorithm \ref{alg:IRsimples} considering $L_1=MSE_f$ and $L_2=MSE_u$. Left: History of functions $L_1$, $L_2$ and $L_1+L_2$. Right: Solution $u(t,x)$ estimated by the Neural Network. Bottom: Solution obtained by the Neural Network vs real solution of Heat equation for specific moments of time. }
        \label{fig_Testes_2_1_custo_campo_curvas}
    \end{figure}

    We observe that the loss functions still show some oscillations in the last epochs. However, they are less abrupt and frequent.

    \item Second set of tests: The maximum number of internal iterations is set $it_{max}=150$, and the parameters $\alpha$ and $\beta$ used are $\alpha=1.5$ and $\beta=1$. Results are shown in Figure \ref{fig_Testes_2_2_custo_campo_curvas}.

    \begin{figure}[!t]
        \centering
        \subfloat[Loss functions.]{\includegraphics[width=1.45in]{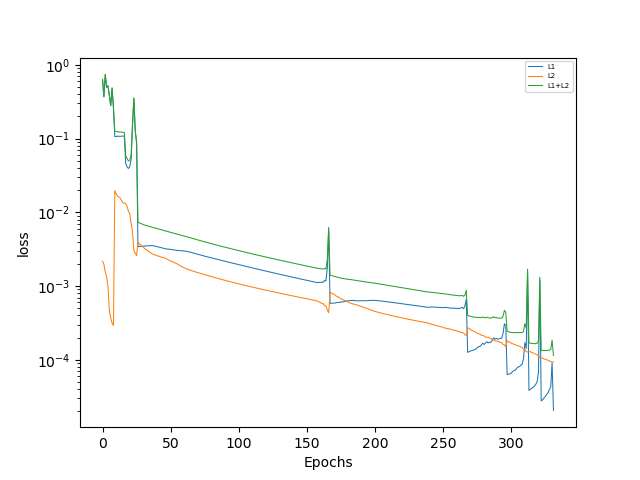}%
        \label{Testes_2_2_custo}}
        \subfloat[Estimated solution $u(t,x)$.]{\includegraphics[width=1.9in]{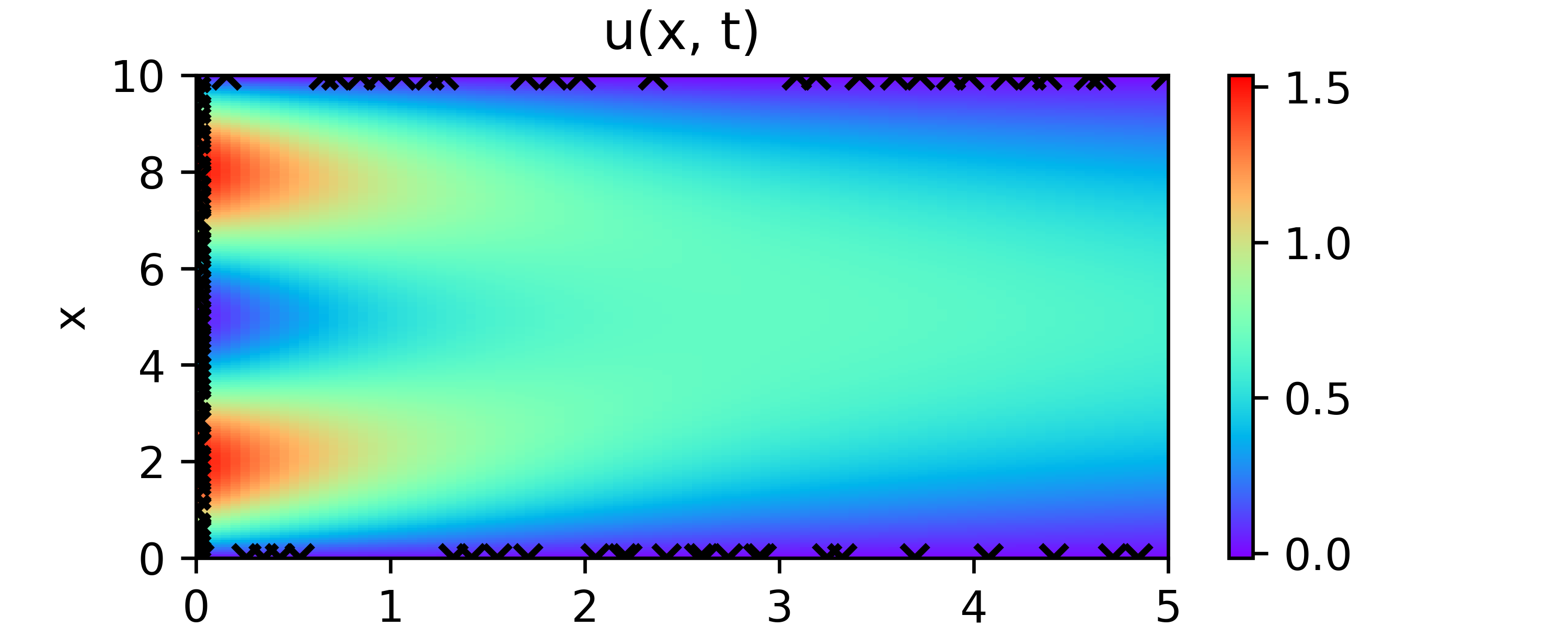}%
        \label{Testes_2_2_campo}}
        
        \subfloat[Estimated solution $u(t,x)$ vs real solution.]{\includegraphics[width=2.9in]{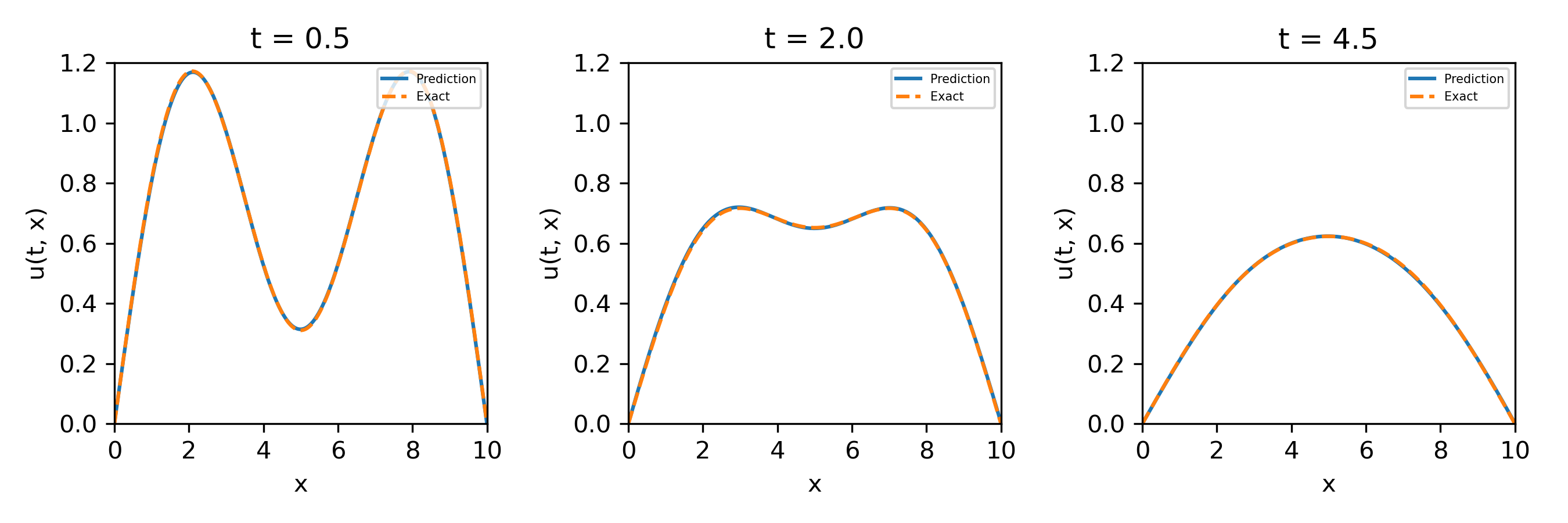}%
        \label{Testes_2_2_curvas}}
        \caption{Results for the Neural Network trained following the Algorithm \ref{alg:IRsimples} considering $L_1=MSE_u$ and $L_2=MSE_f$. Left: History of functions $L_1$, $L_2$ and $L_1+L_2$. Right: Solution $u(t,x)$ estimated by the Neural Network. Bottom: Solution obtained by the Neural Network vs real solution of the Heat equation for specific times.}
        \label{fig_Testes_2_2_custo_campo_curvas}
    \end{figure}
\end{itemize}

\begin{figure}[!t]
\centering
\subfloat[Adam.]{\includegraphics[width=2.9in]{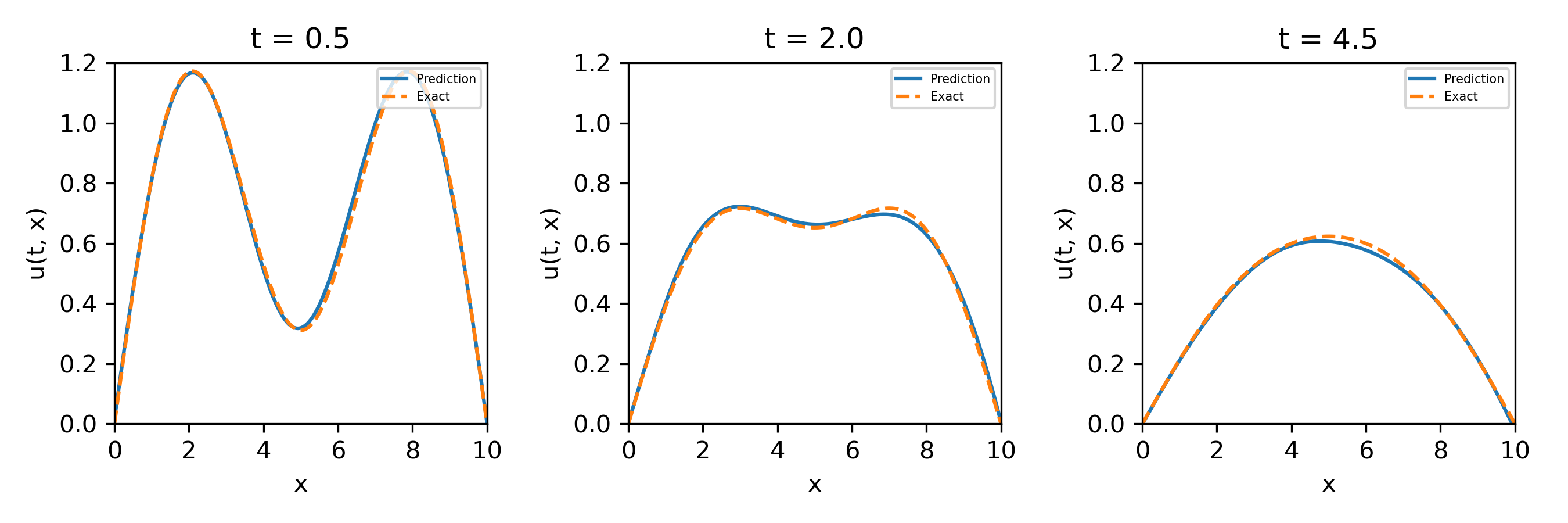}%
\label{comp_1}}

\subfloat[$L_1=MSE_f,L_2=MSE_u$.]{\includegraphics[width=2.9in]{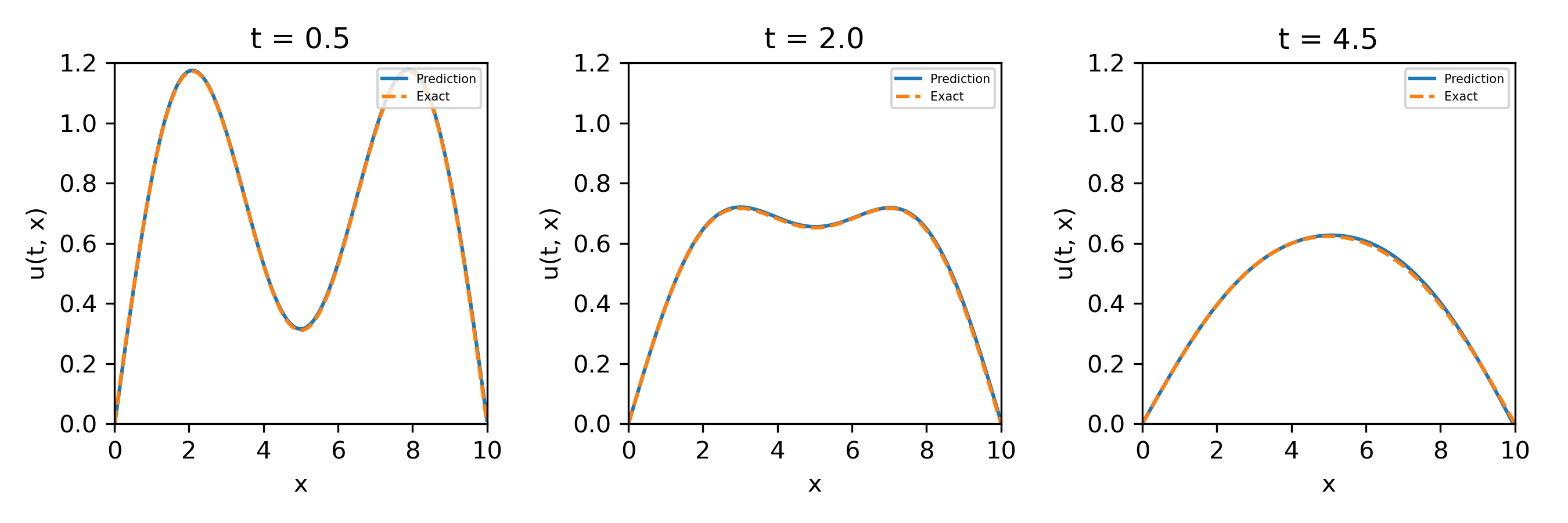}%
\label{comp_2}}

\subfloat[$L_1=MSE_u,L_2=MSE_f$.]{\includegraphics[width=2.9in]{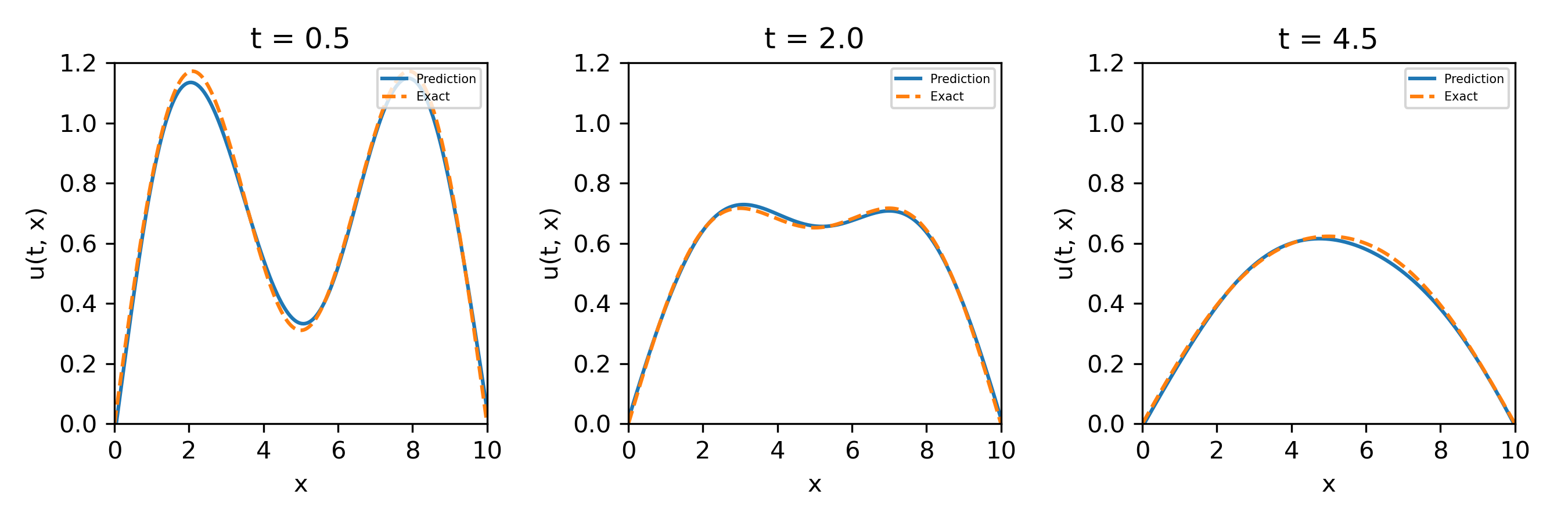}}%
\label{comp_3}
\caption{Solution obtained by the Neural Network vs real solution of  Heat equation for a Neural Network trained with each method after 1000 epochs (or internal iterations).}
\label{comp}
\end{figure}

Finally, to show the speed of convergence, in Figure \ref{comp} we show the solution estimate made by the network for each of the methods after reaching 1000 epochs (or internal iterations).

From here, it is clear that the convergence speed did not suffer significant losses; on the contrary, in Figure \ref{comp_2}, we have a better fit with the solution.

The code used to generate all the examples discussed in this work is available at \url{https://colab.research.google.com/drive/1-8JPbKXQxttt9WtHzYRJtD4kAHcj2Hlf?usp=sharing}.

\section{Conclusions and Future Work}

Algorithm \ref{alg:IRsimples} consistently demonstrated performance at least on par with traditional training methods in all our tests. The imposition of descent conditions reduces the quantity and intensity of oscillations of the loss functions, preventing the loss of solution quality during training. This implies a more stable optimization process, which is one of the proposed method's main advantages and one of its principal contributions. 

We observed that the choice of collocation points significantly influences the behavior of the loss function during training. However, convergence was achieved in all cases. This behavior piques our curiosity and invites further exploration.

More tests are necessary to establish a solid conclusion. However, this work sets a precedent that opens up various avenues for exploration in future research.

The auxiliary optimization methods used within the phases of Algorithm \ref{alg:IRsimples} in this work's tests were chosen due to their similarity to the optimization methods currently used in network training. However, this approach implies the emergence of two extra hyperparameters ($\alpha$ and $\beta$) that must be determined. Besides, there is no guarantee that these methods are ideal. Therefore, we consider Algorithm \ref{alg:IRsimples} promising, and its performance may improve significantly using appropriate methods. Thus, developing appropriate optimization techniques for each phase constitutes a fertile field for future research. In particular, since trust-region methods are widely used and successful in RI algorithms, which inspired this work, trust-region methods must be studied. Future work could also explore methods relevant to Neural Networks or more specific loss functions. 

In the tests presented in this work, the stopping criterion of Algorithm \ref{alg:IRsimples} is the maximum number of iterations. However, in IR algorithms, which inspired this work,  it is expected to use stopping criteria related to the gradient norm of the functions involved or the step size. The study and implementation of more appropriate stopping criteria can be explored in future works.

On the other hand, Algorithm \ref{alg:IRsimples} is not limited to training Physics Informed Neural Networks and can be used in any network whose objective can be represented by two term loss functions. In particular, it could be used in training networks for multitask learning problems. Thus, future research could test the performance of Algorithm \ref{alg:IRsimples} in training networks with other objectives.

\bibliographystyle{IEEEtran}
\bibliography{example}

\end{document}